\documentclass[11pt,preprint]{imsart}
\usepackage{xr}
\usepackage{amsmath,amssymb}
\usepackage{parskip}
\usepackage[numbers]{natbib}
\usepackage{marvosym}
\usepackage[hang,small,bf]{caption}
\usepackage{mathtools,bm}
\usepackage[colorlinks]{hyperref}
\hypersetup{
	colorlinks,%
	citecolor=blue,%
	filecolor=black,%
	linkcolor=blue,%
	urlcolor=blue
}
\usepackage{enumerate}
\usepackage{multirow}
\RequirePackage[OT1]{fontenc}
\usepackage[top=2.5cm, bottom=2.5cm, left=2.5cm, right=2.5cm]{geometry}

\usepackage{amsthm}
\usepackage{amsmath}
\usepackage{amssymb}
\usepackage{thmtools} 
\usepackage{subcaption}
\usepackage{textcomp}

\usepackage{graphicx}
\usepackage{verbatim}
\usepackage{array, float}
\usepackage{fontenc}
\usepackage[toc,page]{appendix}
\usepackage[nottoc]{tocbibind}
\usepackage[english]{babel}
\usepackage{marvosym}

\usepackage[pagewise]{lineno}

\usepackage{mathtools}
\allowdisplaybreaks
\usepackage{bbm}
\usepackage[noabbrev,capitalize]{cleveref}

\newtheoremstyle{exampstyle}
{8pt} 
{8pt} 
{\it} 
{} 
{\bfseries} 
{.} 
{.5em} 
{} 

\theoremstyle{exampstyle}

\newtheorem{theorem}{Theorem}[section]
\newtheorem{proposition}[theorem]{Proposition}
\newtheorem{lemma}[theorem]{Lemma}
\newtheorem{corollary}[theorem]{Corollary}
\newtheorem{remark}{Remark}[section]
\declaretheorem[name={Definition}, sibling=remark]{defn}
\declaretheorem[name={Example}, sibling=theorem]{example}

\newtheoremstyle{remarkstyle}
{8pt} 
{8pt} 
{} 
{} 
{\bfseries} 
{.} 
{.5em} 
{} 

\theoremstyle{remarkstyle}

\usepackage{tikz}
\tikzset{
	treenode/.style = {shape=rectangle, rounded corners,
		draw, align=center,
		top color=white, bottom color=blue!20},
	root/.style     = {treenode, font=\Large, bottom color=yellow},
	env/.style      = {treenode, font=\ttfamily\normalsize},
	con/.style      = {treenode, font=\ttfamily, bottom color=green!25},
	nocon/.style    = {treenode, font=\ttfamily, bottom color=red!30},
	dummy/.style    = {circle,draw}
}

\startlocaldefs

\newcommand{\eat}[1]{}

\DeclareMathOperator*{\argmin}{\arg\!\min}

\renewcommand{\bar}[1]{\overline{#1}}
\renewcommand{\hat}[1]{\widehat{#1}}
\renewcommand{\tilde}[1]{\widetilde{#1}}
\newcommand{\E}{\mathbb{E}}
\renewcommand{\P}{\mathbb{P}}
\newcommand{\Ptac}{\mathcal{P}_2(\R^d)}
\newcommand{\Bmn}{\rho_0}
\newcommand{\Bmnz}{\hat{\rho}_0}

\newcommand{\Pac}{\mathcal{P}_{\mathrm{ac}}(\R^d)}
\newcommand{\Paco}{\mathcal{P}_{\mathrm{ac}}(\R^{d_1})}
\newcommand{\Pact}{\mathcal{P}_{\mathrm{ac}}(\R^{d_2})}

\newcommand{\mx}{x}
\newcommand{\my}{y}
\newcommand{\ttg}{\tilde{T}_{m,n}^{\gamma}}
\newcommand{\Tes}{\tilde{T}_{m,n}}

\newcommand{\R}{\mathbb{R}}

\newcommand{\pn}{\pi^{(n)}}
\newcommand{\mn}{\mu^{(n)}}
\newcommand{\nn}{\nu^{(n)}}

\newcommand{\Hy}{\mathrm{H}}

\newcommand{\B}{\mathcal{B}}
\newcommand{\cS}{\mathcal{S}}
\newcommand{\phmn}{\varphi_0}
\newcommand{\pst}{\Psi_{\tmu,\tnu}}
\newcommand{\psw}{\Psi_{\tmu,\bar{\nu}_m}}

\newcommand{\fm}{f_{\mu}}
\newcommand{\fn}{f_{\nu}}
\newcommand{\bph}{\boldsymbol{\Phi}}
\newcommand{\bps}{\boldsymbol{\Psi}}

\newcommand{\bX}{X}
\newcommand{\bY}{Y}

\newcommand{\tgm}{\tilde{\Gamma}_{\mathrm{min}}}

\newcommand{\xmo}{\mathcal{X}_{n,\mu}}

\newcommand{\xn}{\mathcal{Y}}

\newcommand{\mm}{\mathfrak{m}}

\newcommand{\tM}{\hat{\mu}_{m,M}}
\newcommand{\tN}{\hat{\nu}_{n,M}}

\newcommand{\psm}{\Psi_{\mu,\nu}}

\newcommand{\emu}{\tilde \mu_m}
\newcommand{\enu}{\tilde \nu_n}
\newcommand{\hmu}{\hat \mu_m}
\newcommand{\hnu}{\hat \nu_n}
\newcommand{\tmu}{\tilde \mu_m}
\newcommand{\tnu}{\tilde \nu_n}

\newcommand{\hfm}{\hat{f}_{\mu}}
\newcommand{\hfn}{\hat{f}_{\nu}}

\newcommand{\Tmn}{T_0}
\newcommand{\eps}{\varepsilon}
\newcommand{\hao}{\hat{A}_{n,1}}
\newcommand{\had}{\hat{A}_{n,2}}
\newcommand{\har}{\hat{A}_{n,3}}
\newcommand{\harr}{\hat{A}_{n,3}^{\mathrm{OR}}}
\newcommand{\hadr}{\hat{A}_{n,2}^{\mathrm{OR}}}
\newcommand{\haor}{\hat{A}_{n,1}^{\mathrm{OR}}}

\newcommand{\bmor}{\Bmn^{\mathrm{OR}}}
\newcommand{\G}{\mathcal{G}}

\newcommand{\cN}{\mathcal{N}}

\usepackage{hyperref}

\newcommand{\tfm}{\tilde{f}_{\mu}}
\newcommand{\tfn}{\tilde{f}_{\nu}}

\newcommand{\hmn}{\hat{\mu}_n}
\newcommand{\ohs}{\hhs^{\mathrm{OR}}}
\newcommand{\hxr}{\hat{x}^{\mathrm{OR}}}
\newcommand{\hyr}{\hat{y}^{\mathrm{OR}}}

\newcommand{\U}{\mathrm{Unif}}

\newcommand{\X}{\mathcal{X}}
\newcommand{\F}{\mathcal{F}}

\newcommand{\Y}{\mathcal{Y}}

\newcommand{\hu}{\hat{u}_n}
\newcommand{\hv}{\hat{v}_n}

\newcommand{\oto}{T_1^{(n)}}
\newcommand{\ott}{T_2^{(n)}}

\newcommand{\hto}{\hat{T}_{1,n}}
\newcommand{\htt}{\hat{T}_{2,n}}
\newcommand{\hhs}{\hat{\mathrm{rHSIC}}}
\newcommand{\nhs}{\mathrm{rHSIC}}

\newcommand{\hx}{\hat{x}}
\newcommand{\hy}{\hat{y}}

\theoremstyle{plain}

\numberwithin{equation}{section}
\newtheorem{Assumption}{Assumption}

\def\beq{\begin{equation}}
\def\eeq{\end{equation}}
\def\ba{\begin{enumerate}[(a)]}
	\def\bei{\begin{enumerate}[(i)]}
		\def\be{\begin{enumerate}[(1)]}
			\def\ee{\end{enumerate}}
		\def\bi{\begin{itemize}}
			\def\ei{\end{itemize}}
		\def\beg{\begin{eg}}
			\def\eeg{\end{eg}}
		\def\bd{\begin{defn}}
			\def\ed{\end{defn}}
		\def\bt{\begin{thm}}
			\def\et{\end{thm}}
		\def\bl{\begin{lemma}}
			\def\el{\end{lemma}}
		\def\bfac{\begin{fact}}
			\def\efac{\end{fact}}
		
		\def\bc{\begin{corollary}}
			\def\ec{\end{corollary}}
		\def\bp{\begin{prop}}
			\def\ep{\end{prop}}
		\def\bo{\begin{observe}}
			\def\eo{\end{observe}}
		\def\bas{\begin{assumption}}
			\def\eas{\end{assumption}}

\endlocaldefs

\usepackage{abstract}
\usepackage{xr}
\usepackage[numbers]{natbib}
\begin{document}
	
	\renewcommand{\abstractname}{}    
	\renewcommand{\absnamepos}{empty}
	
	\begin{frontmatter}
		\title{Rates of Estimation of Optimal Transport Maps using Plug-in Estimators via  Barycentric Projections}
		
		\runtitle{Rates of Convergence of Barycentric Projections}		
		
		\begin{aug}
			\author{\fnms{Nabarun} 
				\snm{Deb}, \ead[label=e1]{nd2560@columbia.edu}}\author{\fnms{Promit} \snm{Ghosal,}
				\ead[label=e2]{promit@mit.edu}}
			\and
			\author{\fnms{Bodhisattva} \snm{Sen}\thanksref{t3}\ead[label=e3]{bodhi@stat.columbia.edu}}
			\affiliation{
				Columbia University\thanksmark{a1},
				Massachusetts Institute of Technology\thanksmark{a2}, and
				Columbia University\thanksmark{a3}
			}
			\thankstext{t3}{Supported by NSF grant DMS-2015376.}
			
			\runauthor{Deb, Ghosal, and Sen}
			
			\address{1255 Amsterdam Avenue \\
				New York, NY 10027\\
				\printead{e1} \\
			}
			\address{77 Massachusetts Ave \\
				Cambridge, MA 02139 \\
				\printead{e2} \\			
			}
			\address{1255 Amsterdam Avenue \\
				New York, NY 10027\\
				\printead{e3} \\
			}
		\end{aug}
		\vspace{0.2in}
		\begin{abstract}
			    Optimal transport maps between two probability distributions $\mu$ and $\nu$ on $\R^d$ have found extensive applications in both machine learning and statistics. In practice, these maps need to be estimated from data sampled according to $\mu$ and $\nu$. Plug-in estimators are perhaps most popular in estimating transport  maps in the field of computational optimal transport. In this paper, we provide a comprehensive analysis of the rates of convergences for general plug-in estimators defined via barycentric projections. Our main contribution is a new stability estimate for barycentric projections which proceeds under minimal smoothness assumptions and can be used to analyze general plug-in estimators. We illustrate the usefulness of this stability estimate by first providing rates of convergence for the natural discrete-discrete and semi-discrete estimators of  optimal transport maps. We then use the same stability estimate to show that, under additional smoothness assumptions of Besov type or Sobolev type, wavelet based or kernel smoothed plug-in estimators respectively speed up the rates of convergence and significantly mitigate the curse of dimensionality suffered by the natural discrete-discrete/semi-discrete estimators. As a by-product of our analysis, we also obtain faster rates of convergence for plug-in estimators of $W_2(\mu,\nu)$, the Wasserstein distance between $\mu$ and $\nu$, under the aforementioned smoothness assumptions, thereby complementing recent results in Chizat et al. (2020). Finally, we illustrate the applicability of our results in obtaining rates of convergence for Wasserstein barycenters between two probability distributions and obtaining asymptotic detection thresholds for some recent optimal-transport based tests of independence.  
		\end{abstract}

		
	\end{frontmatter}

	\section{Introduction}
	Given two random variables $\bX\sim\mu$ and $\bY\sim\nu$, where $\mu,\nu$ are probability measures on $\R^d$, $d\geq 1$, the problem of finding a ``nice" map $\Tmn(\cdot)$ such that $\Tmn(\bX)\sim\nu$ has numerous applications in machine learning such as domain adaptation and data integration~\cite{gopalan2011domain,forrow2019statistical,courty2016optimal,courty2017,damodaran2018deepjdot,seguy2018large}, dimension reduction~\cite{gunsilius2019independent,bigot2017geodesic,masarotto2019procrustes}, generative models~\cite{goodfellow2014generative,kingma2013auto,li2015generative,salimans2018improving}, to name a few. Of particular interest is the case when $\Tmn(\cdot)$ is obtained by minimizing a cost function, a line of work initiated by Gaspard Monge~\cite{Monge1781} in 1781 (see~\eqref{eq:OTM} below), in which case $\Tmn(\cdot)$ is termed an \emph{optimal transport} (OT) map and has applications in shape matching/transfer problems~\cite{ferradans2014regularized,su2015optimal,chizat2018scaling,reich2011dynamical}, Bayesian statistics~\cite{reich2013nonparametric,el2012bayesian,Kandaswamy2018,kim2013efficient}, econometrics~\cite{galichon2016optimal,blanchet2016optimal,chiappori2010hedonic,friesz1979model,ekeland2010optimal}, nonparametric statistical inference~\cite{deb2021multivariate,shi2020distribution,Shi2020,deb2020measuring,deb2021efficiency}; also see~\cite{V03,V09,santambrogio2015optimal} for book-length treatments on the subject. In this paper, we will focus on the OT map obtained using the standard Euclidean cost function, i.e.,
	\begin{equation}\label{eq:OTM}
	\Tmn:=\argmin\limits_{T:T\#\mu=\nu}\E\lVert \bX-T(\bX)\rVert^2,
	\end{equation}
	where $T\#\mu=\nu$ means $T(\bX)\sim\nu$ for $\bX\sim\mu$. The estimation of $\Tmn$ has attracted a lot of interest in recent years due to its myriad applications (as stated above) and interesting geometrical properties (see~\cite{McCann95,Gigli2011,Brenier1991} and~\cref{def:otmpm} below). In practice, the main hurdle in constructing estimators for  $\Tmn$ is that the explicit forms of the measures $\mu,\nu$ are unknown; instead only random samples $$X_1,\ldots ,X_m\sim \mu \qquad\quad\mbox{and} \qquad\quad Y_1,\ldots ,Y_n\sim\nu$$ are available. A natural strategy in this scenario is to estimate $\Tmn$ using $\Tes$, where $\Tes$ is computed as in~\eqref{eq:OTM} with $\mu$ and  $\nu$ replaced by $\emu$ and $\enu$ which are empirical approximations of $\mu$ and $\nu$ based on $X_1,\ldots ,X_m$ and $Y_1,\ldots ,Y_n$ respectively (see~\cref{def:bcenterproj}). Such estimators are often called \emph{plug-in estimators} and have been used extensively, see~\cite{Sommerfeld2018,Merigot11,Merigot2020,benamou2000computational,papadakis2014optimal,gunsilius2018convergence,chizat2020faster}. 
	
	The main goal of this paper is to study the \emph{rates of convergence} of general plug-in estimators of $\Tmn$ under a unified framework. We show that when $\emu$ and $\enu$ are chosen as $\hmu$ and $\hnu$ respectively, where $\hmu$ and $\hnu$ are the standard empirical distributions supported on $m$ and $n$ atoms, i.e., \begin{equation}\label{eq:empdef}\hmu:=\frac{1}{m}\sum_{i=1}^m \delta_{X_i}\qquad\quad \mbox{and} \qquad\quad \hnu:=\frac{1}{n}\sum_{j=1}^n \delta_{Y_j},\end{equation}
	$\Tes$ (appropriately defined using~\cref{def:bcenterproj}) converges at a rate of $m^{-2/d}+n^{-2/d}$ for $d\geq 4$. This rate happens to be minimax optimal under minimal smoothness assumptions (see~\cite[Theorem 6]{hutter2021minimax}) but suffers from the \emph{curse of dimensionality}. We next show that, if $\mu$ and $\nu$ are known to admit sufficiently smooth densities, it is possible to apply wavelet or kernel based smoothing techniques on $\hmu$ and $\hnu$ to obtain plug-in estimators that mitigate the aforementioned curse of dimensionality. 
	
	Our next contribution pertains to the estimation of $W_2^2(\mu,\nu)$ (the squared Wasserstein distance), see~\eqref{eq:wass} below, a quantity of independent interest in statistics and machine learning with applications in structured prediction~\cite{Frogner2015,Luise2018}, image analysis~\cite{glaunes2004diffeomorphic,bonneel2011displacement}, nonparametric testing~\cite{boeckel2018multivariate,Ramdas2017wasserstein}, generative modeling~\cite{shakir2017,bernton2017inference}, etc. In this paper, we also obtain rates of convergence for plug-in estimators $W_2^2(\tmu,\tnu)$ of $W_2^2(\mu,\nu)$. We show that kernel smoothing $\hmu$ and $\hnu$ can be used to obtain plug-in estimators of $W_2^2(\mu,\nu)$ that mitigate the curse of dimensionality as opposed to a direct plug-in approach using $\hmu$ and $\hnu$ (as used in~\cite[Theorem 2]{chizat2020faster}). This provides an answer to the open question of estimating $W_2^2(\mu,\nu)$ when $\mu$, $\nu$ admit smooth densities laid out in~\cite{chizat2020faster}. 
	\subsection{Background on optimal transport}\label{sec:setting}
	In this section, we present some basic concepts and results associated with the OT problem that will play a crucial role in the sequel. Let $\Pac$ denote the set of all Lebesgue absolutely continuous probability measures on $\R^d$ and $\Ptac$ be the set of probability measures with finite second moments.  Then the $2$-\emph{Wasserstein} distance (squared) between $\mu,\nu\in\Ptac$ is defined as:
	\begin{equation}\label{eq:wass}
	W_2^2(\mu,\nu):=\min\limits_{\pi\in \Pi(\mu,\nu)}\int \lVert \mx-\my\rVert^2\,d\pi(\mx,\my),
	\end{equation} 
	where $\Pi(\mu,\nu)$ is the set of probability measures on $\R^d\times \R^d$ with marginals $\mu$ and $\nu$. The optimization problem in~\eqref{eq:wass} is often called the \emph{Kantorovich relaxation} (see~\cite{Kantorovic1948,Kantorovic1942}) of the optimization problem in~\eqref{eq:OTM}. The existence of a minimizer in~\eqref{eq:wass} follows from~\cite[Theorem 4.1]{V09}.
	
	\begin{proposition}[Brenier-McCann polar factorization theorem,~see~\cite{V03,McCann95}]\label{prop:bmopt}
		Suppose $\mu,\nu\in\Pac$. Then there exists a $\mu$-a.e.~(almost everywhere) unique function $\Tmn(\cdot):\R^d\to\R^d$, which is the gradient of a real-valued $d$-variate convex function, say $\phmn(\cdot):\R^d\to\R$, such that $\Tmn\#\mu=\nu$. Further, the distribution defined as $\pi(A\times B)=\mu(A\cap (\Tmn)^{-1}(B))$ for all Borel sets $A,B\subseteq\R^d$ is the unique minimizer in~\eqref{eq:wass} provided $\mu,\nu\in\Pac\cap\Ptac$.
	\end{proposition}
	\begin{defn}[OT map and potential function]\label{def:otmpm}
		The function $\Tmn:\R^d\to\R^d$ in~\cref{prop:bmopt} which satisfies $\Tmn\#\mu=\nu$ will be called the  OT map from $\mu$ to $\nu$. The convex function $\phmn(\cdot)$ in~\cref{prop:bmopt} satisfying $\nabla \phmn=\Tmn$ will be termed the OT potential.
	\end{defn}
	
	The next and final important ingredient is the alternate dual representation of~\eqref{eq:wass} which gives:
	\begin{align}
	\frac{1}{2}W_2^2(\mu,\nu)&=\frac{1}{2}\int \lVert x\rVert^2\,d\mu(x)+\frac{1}{2}\int \lVert y\rVert^2\,d\nu(y)-\min_{f\in\F} \cS_{\mu,\nu}(f),\,\qquad \mbox{where}\label{eq:altwass} \\ \cS_{\mu,\nu}(f)&=\int f\,d\mu+\int f^*\,d\nu.\label{eq:dualf}
	\end{align}
	Here $\F$ denotes the space of convex functions on $\R^d$ which are also elements of $L^1(\mu)$ and $f^*(\cdot)$ is the standard Legendre-Fenchel dual defined as:
	\begin{equation}\label{eq:lfdual}
	f^*(x) := \sup_{y \in \R^d} [y^\top x - f(y)], \qquad \mbox{for } x \in \mbox{dom}(f).
	\end{equation}
	
	
	
	\subsection{Estimating OT map via barycentric projection}\label{sec:baryproj}
	Recall the setting from the Introduction. Let $\tmu,\tnu\in\Ptac$. Here $\tmu,\tnu$ need not be absolutely continuous and can be very general. Intuitively, $\tmu$ and $\tnu$ should be viewed as some empirical approximation of $\mu$ and $\nu$ respectively. 
	
	\begin{example}[Simple choices of $\tmu$ and $\tnu$]\label{ex:baseplug}
		Let $\bX_1,\ldots ,\bX_m\overset{i.i.d.}{\sim}\mu$ and $\bY_1,\ldots ,\bY_n\overset{i.i.d.}{\sim}\nu$; in which case a natural choice would be to set $\tmu=\hmu$ and $\tnu=\hnu$ where $\hmu$ and $\hnu$ are the empirical distributions on $\bX_1,\ldots ,\bX_m$ and $\bY_1,\ldots ,\bY_n$ respectively, as defined in~\eqref{eq:empdef}. This is the standard choice adopted in the discrete-discrete Kantorovich relaxation; see~\cite[Section 2.3]{peyre2019computational}. Another popular choice is $\tmu=\hmu$, $\tnu=\nu$ or $\tmu=\mu$, $\tnu=\hnu$. This is the semi-discrete Kantorovich problem and is popular when one of the measures is fully  specified; see~\cite{Chernozhukhov2017,ghosal2019multivariate}.
	\end{example}
	
	A natural way to estimate $\Tmn(\cdot)$, as defined in~\eqref{eq:OTM}, would be to approximate it using the OT map from $\tmu$ to $\tnu$. However as $\tmu$ and $\tnu$ may not be elements of $\Pac$,~\cref{prop:bmopt} does not apply and an OT map \emph{may not exist} from $\tmu$ to $\tnu$. Such is the case in~\cref{ex:baseplug} in the \emph{discrete-discrete} case when $m\neq n$. To circumvent this issue, we leverage the notion of barycentric projections (see~\cite[Definition 5.4.2]{Ambro08}) defined below:
	\begin{defn}[Barycentric projection]\label{def:bcenterproj}
		Define the set
		$$\tgm:=\argmin\limits_{\pi\in\Pi(\tmu,\tnu)} \int \lVert \mx-\my\rVert^2\,d\pi(\mx,\my).$$
		The optimization problem above is the plug-in analog of the optimization problem on the right hand side of~\eqref{eq:wass}. Given any $\gamma\in\tgm$, define the barycentric projection of $\gamma$ as the conditional mean of $\my$ given $\mx$ under $\gamma$, i.e.,
		\begin{equation}\label{eq:conmean}\Tes(x)\equiv \ttg(\mx):=\frac{\int_{\my} \my \,d\gamma(\mx,\my)}{\int_{\my} d\gamma(\mx,\my)},\qquad \mbox{for}\ x\in \mbox{supp}\left({\tmu}\right).\end{equation}
	\end{defn}
	In general, $\tgm$ need not be a singleton which is why we index the barycentric projection $\ttg(\cdot)$ by $\gamma\in\tgm$. Note that $\ttg(\cdot)$ need not be a transport map; however, if an OT map exists then it must be equal to $\ttg(\cdot)$ ($\tmu$-a.e.). Our goal is to obtain stochastic upper bounds for
	\begin{equation}\label{eq:OTrate}
	\sup\limits_{\gamma\in\tgm} \int \big\lVert \ttg(\mx)-\Tmn(\mx)\big\rVert^2\,d\tmu(\mx).
	\end{equation}
	In addition, our proof techniques also yield rates of convergence for
	\begin{equation}\label{eq:wassrate}
	\big|W_2^2(\tmu,\tnu)-W_2^2(\mu,\nu)\big|.
	\end{equation}
	In this paper, we will focus on $d\geq 2$. Due to the canonical ordering of $\R$, the case $d=1$ can be handled easily using the classical Hungarian embedding theorem~\cite{Komlos1976}.
	\subsection{Contributions}\label{sec:contrib}
	\begin{enumerate}
		\item We provide a new and flexible stability estimate~\cref{thm:newubd} which yields a unified approach to obtaining \emph{rates} of convergence for general plug-in estimators of the OT map $\Tmn(\cdot)$. Unlike existing stability estimates,~\cref{thm:newubd} holds for the barycentric projection (which is the same as the OT map when it exists) and does not require any smoothness assumptions on $\tmu$, $\tnu$ or $\ttg(\cdot)$; also see~\cref{rem:comstab} for a comparison with the existing literature.
		
		\item in Sections~\ref{sec:nsmooth}~and~\ref{sec:smooth}, we use~\cref{thm:newubd}  to bound~\eqref{eq:OTrate}~and~\eqref{eq:wassrate}:
		\begin{itemize}
			\item In~\cref{sec:nsmooth}, we show that in both the discrete-discrete and semi-discrete Kantorovich relaxation problems (see~\cref{ex:baseplug}), the rate of convergence of~\eqref{eq:OTrate} is $m^{-2/d}+n^{-2/d}$ for $d\geq 4$ when $\Tmn$ is assumed to be Lipschitz (see~\cref{thm:nsmooth}), which is the minimax rate (see~\cite[Theorem 6]{hutter2021minimax}). To the best of our knowledge, rates of convergence for these natural estimators weren't previously established in the literature. 
			
			\item In~\cref{sec:smooth}, we show that the curse of dimensionality in the above rates can be mitigated provided $\mu$ and $\nu$ admit Besov smooth densities (see~\cref{sec:wavelet}) or (uniform) Sobolev smooth densities (see~\cref{sec:kernel}). In~\cref{sec:wavelet}, our plug-in estimator is obtained using natural wavelet based density estimators. The rate of convergence in~\eqref{eq:OTrate} turns out to be $n^{-\frac{1+s}{d+2s}}$ where $s$ denotes the degree of Besov smoothness (see~\cref{thm:smoothwav}). Note that by choosing $s$ large enough, the exponent in the rate can be made arbitrarily close to $1/2$, thereby reducing the curse of dimensionality. In~\cref{sec:kernel}, our plug-in estimator is obtained by choosing $\tmu$ (and $\tnu$) as the convolution of $\hmu$ (and $\hnu$) and a smooth kernel with an appropriate bandwidth. Under this choice, the rate of convergence in~\eqref{eq:OTrate} is $m^{-\left(\frac{s+2}{d}\wedge \frac{1}{2}\right)}+n^{-\left(\frac{s+2}{d}\wedge \frac{1}{2}\right)}$, where $s$ denotes the degree of Sobolev smoothness (see~\cref{thm:smooth}). Clearly, if $2(s+2)\geq d$, the rate of convergence becomes dimension-free and mitigates the curse of dimensionality. We also show the same rates of convergence mentioned above also hold for~\eqref{eq:wassrate} (see e.g.,~\cref{prop:wassrate}) which makes a strong case in favor of incorporating smoothness in the construction of plug-in estimators as was conjectured in~\cite{chizat2020faster}.
		\end{itemize}
		
		\item In~\cref{sec:dismooth}, we use a discretization technique from~\cite{weed2019estimation} to construct discrete approximations to the smoothed $\tmu$ and $\tnu$ from the previous paragraph that in turn yield computable plug-in estimators for $\Tmn$ (provided one can sample from $\tmu$ and $\tnu$) that also achieve the same statistical guarantees as the smoothed plug-in estimator from~\cref{sec:smooth} (see~\cref{thm:dismooth}). However the number of atoms required in the discretizations and correspondingly the computational complexity increases with the degree of smoothness; this highlights a statistical and computational trade off.
		
		\item We provide implications of our results in popular applications of OT such as estimating the barycenter of two multivariate probability distributions (see~\cref{thm:barate} in~\cref{sec:bar}) and in nonparametric independence testing (see~\cref{thm:dethresh} in~\cref{sec:nonptest}).
	\end{enumerate}
	\subsection{Related work}\label{sec:relwork}
	
	Many recent works have focused on obtaining consistent estimators of $\Tmn$ using the plug-in principle, see~\cite{Chernozhukhov2017,ghosal2019multivariate} (in the semi-discrete problem) and~\cite{hallin2021distribution,zemel2019frechet,deb2021multivariate} (in the discrete-discrete problem). In~\cite{ghosal2019multivariate}, the authors have studied the rate of convergence of the semi-discrete optimal transport map from $\nu$ (absolutely continuous) to $\hmu$. This paper complements the aforementioned  papers by studying the rates of convergence for general plug-in estimators in a unified fashion. In two other papers~\cite[Theorem 1.1]{berman2020convergence}~and~\cite[Section 4]{li2020quantitative}, the authors use a ``Voronoi tessellation" approach to estimate $\Tmn$, however the rates obtained in this paper, even in the absence of smoothness, are strictly better than those in~\cite{berman2020convergence,li2020quantitative}. Perhaps the most closely related paper to ours would be~\cite{gunsilius2018convergence}. In~\cite{gunsilius2018convergence}, the author uses variational techniques to arrive at stability estimates while we exploit the Lipschitz nature of the OT map (see~\cref{def:otmpm}). Further the rates in this paper have exponents  $\frac{s+2}{d}\wedge \frac{1}{2}$ which are \emph{strictly better} than the exponents $\frac{s+2}{2(s+2)+d}$ obtained in~\cite[Proposition 1]{gunsilius2018convergence} under the same smoothness assumptions  (Sobolev type of order $s$, see~\cref{def:holder}). In another line of work~\cite{hutter2021minimax}, the authors use theoretical wavelet based estimators (not of the plug-in type) of $\Tmn$ to obtain nearly minimax optimal rates of convergence. However these estimators, by themselves, are not transport maps between two probability measures, which makes them harder to interpret. In contrast, our focus is on obtaining rates of convergence for plug-in estimators, which are transport maps between natural aprroximations of $\mu$ and $\nu$. Such plug-in type strategies are a lot more popular in computational OT~\cite{Sommerfeld2018,Merigot11,Merigot2020,benamou2000computational,papadakis2014optimal,gunsilius2018convergence,chizat2020faster}. 
	
	In terms of obtaining rates of convergence for~\eqref{eq:wassrate}, some attempts include~\cite{Sommerfeld2018,Rippl2016} where parametric rates are obtained when $\mu,\nu$ are known to be finitely supported or are both Gaussian. In a related problem, bounds for $W_2^2(\hmu,\mu)$ were obtained in~\cite{Talagrand1994,Barthe2013,Dereich2013,Fournier2015,niles2019estimation,weed2019estimation}. Using these bounds, it is easy to get a $n^{-1/d}$ rate of convergence for~\eqref{eq:wassrate}. This rate was recently improved to $n^{-2/d}$ in~\cite{chizat2020faster} under no smoothness assumptions. Our rates coincide with the $n^{-2/d}$ rate from~\cite{chizat2020faster} under no smoothness assumptions. But further, we show in this paper that the curse of dimensionality in the above rate can be mitigated by incorporating smoothness into the plug-in procedure. 
	\section{Main results}\label{sec:mainres}
	
	Recall the definition of $\phmn(\cdot)$ from~\cref{def:otmpm}. The following is our main result.
	\begin{theorem}[Stability estimate]\label{thm:newubd}
		Suppose that $\mu,\nu\in\Pac\cap\Ptac$ and $\tmu,\tnu\in\Ptac$. Assume that $T_0(\cdot)$ (as defined in~\eqref{eq:OTM}) is $L$-Lipschitz ($L >0$). Then, 
		\begin{align}\label{eq:globaldiffred}
		\sup\limits_{\gamma\in\tgm} \int \lVert \ttg(\mx)-\Tmn(\mx)\rVert^2& \,d\tmu(\mx)\le L \max\left\{\bigg|\int \pst^* \,d(\tnu-\bar{\nu}_m)\bigg|,\bigg|\int \psw^* \,d(\tnu-\bar{\nu}_m)\bigg|\right\}\nonumber \\ &+2L\int \phmn^*(\my)\,d(\tnu-\bar{\nu}_m)(\my),
		\end{align}
		where $\bar{\nu}_m=\Tmn\#\tmu$, $\phmn^*(\cdot)$ is defined as in~\eqref{eq:lfdual}, and with $\cS_{\cdot,\cdot}(\cdot)$ defined as in~\eqref{eq:dualf}, $\pst(\cdot):=\argmin_{f\in\F} \cS_{\tmu,\tnu}(f),\quad \psw(\cdot):=\argmin_{f\in\F} \cS_{\tmu,\bar{\nu}_m}(f)$.
	\end{theorem}
	The proof of~\cref{thm:newubd} (see~\cref{sec:mainrespf}) starts along the same lines as the proof of the curvature estimate in~\cite[Proposition 3.3]{Gigli2011}. This is followed by some careful manipulations of $W_2^2(\cdot,\cdot)$ (as in~\eqref{eq:wass}) and an application of the conditional version of Jensen's inequality, see~\eqref{eq:newubd3}. The final step of the proof uses the dual representation in~\eqref{eq:altwass} with techniques similar to some intermediate steps in the proof of~\cite[Proposition 2]{Gonzalo2019}~and~\cite[Lemma 3]{chizat2020faster}.
	
	\begin{remark}[Comparison with other stability estimates]\label{rem:comstab}
		\cref{thm:newubd} provides some important advantages to existing stability estimates in the literature. One of the earliest results in this direction can be found in~\cite[Proposition 3.3]{Gigli2011} but their bound involves a push-forward constraint which makes it hard to use for rate of convergence analysis. A bound similar to~\cref{thm:newubd} is presented in~\cite[Lemma 5.1]{ghosal2019multivariate} but there the authors assume the existence of an OT map from $\tmu$ to $\tnu$. Therefore, it does not apply to the discrete-discrete problem where $\tmu=\hmu$ and $\tnu=\hnu$ with $m\neq n$. Overcoming all these limitations is an important contribution of~\cref{thm:newubd} and allows us to deal with popular  plug-in estimators all in one go. The stability estimate in~\cite[Proposition 10]{hutter2021minimax} on the other hand requires $\tmu$, $\tnu$ to be sufficiently smooth and hence it does not hold for discrete-discrete or semi-discrete plug-in estimators (see~\cref{ex:baseplug}). Further their result requires all the measures involved to be compactly supported unlike the much milder requirements of~\cref{thm:newubd}. However, a shortcoming of~\cref{thm:newubd} is that it is hard to obtain rates faster than $n^{-1/2}$ using it directly, whereas~\cite{hutter2021minimax} can obtain rates arbitrarily close to $n^{-1}$. This is a price we pay for analyzing natural and popular plug-in estimators as opposed to the (more intractable) wavelet based estimators in~\cite{hutter2021minimax}.
	\end{remark}
	
	\begin{remark}[\emph{How to use~\cref{thm:newubd} to obtain rates of convergence?}]\label{rem:pftech} Note that the second term on the right hand side of~\eqref{eq:globaldiffred}, under appropriate moment assumptions, is $O_p(m^{-1/2}+n^{-1/2})$ (free of dimension) by a direct application of Markov's inequality. We therefore focus on the first term. By~\eqref{eq:dualf}, $\pst^*(\cdot)$, $\psw^*(\cdot)\in \F$. Further, by Caffarelli's regularity theory~\cite{Caf1992,Caf92,Caf1996}, depending on the ``smoothness" of $\tmu$, $\tnu$, it can be shown that there exists a further class of functions $\F_s$ (see Remarks~\ref{rem:pfnstech}~and~\ref{rem:pfstech}) such that $\pst^*(\cdot)$, $\psw^*(\cdot)\in \F\cap \F_s$. Thus, we can bound the first term on the right hand side of~\eqref{eq:globaldiffred} as:
		\begin{small}
			\begin{align}\label{eq:strat}
			\max\left\{\bigg|\int \pst^* \,d(\tnu-\bar{\nu}_m)\bigg|,\bigg|\int \psw^* \,d(\tnu-\bar{\nu}_m)\bigg|\right\}\leq \sup_{f\in\F\cap \F_s}\bigg|\int f\,d(\tnu-\bar{\nu}_m)\bigg|.
			\end{align}
		\end{small}
		The right hand side of~\eqref{eq:strat} can now be bounded using the corresponding  Dudley's entropy integral bounds using empirical process techniques, see~\cite[Lemmas 19.35-19.37]{vdV98}. 
	\end{remark}
	To conclude, the two main steps in our strategy are identifying the family of functions $\F_s$ and computing Dudley's entropy integral. Further, the more the smoothness of $\tmu$, $\tnu$, the smaller is the class of functions $\F_s$ and smaller the supremum on the right hand side of~\eqref{eq:strat}. This shows why better rates can be expected under smoothness assumptions.
	
	\subsection{Natural non-smooth plug-in estimator}\label{sec:nsmooth}
	In this case, we discuss the rates of convergence for the discrete-discrete problem and the semi-discrete problem, where \emph{no smoothness} is available on $\tmu$ and $\tnu$.
	
	\begin{theorem}\label{thm:nsmooth} Suppose that $T_0(\cdot)$ is $L$-Lipschitz,  $\nu$ is compactly supported and $\E \exp(t\lVert X_1\rVert^{\alpha})<\infty$ for some $t>0$, $\alpha>0$. 	
		
		\textbf{(Discrete-discrete):} Set $\tmu=\hmu$ and $\tnu=\hnu$. Then the following holds:
		\begin{equation}\label{eq:nsym}\sup\limits_{\gamma\in\tgm} \int \lVert \ttg(\mx)-\Tmn(\mx)\rVert^2\,d\tmu(x)=O_p\left(r^{(m,n)}_d\times (\log{(1+\max\{m,n\})})^{t_{d,\alpha}}\right),\end{equation}  \begin{equation}\label{eq:nsmrate}\mbox{where}\quad r^{(m,n)}_d:=\begin{cases} m^{-1/2}+n^{-1/2} & \mbox{for}\ d=2,3,\\ m^{-1/2}\log{(1+m)}+n^{-1/2}\log{(1+n)} & \mbox{for}\ d=4,\\ m^{-2/d}+n^{-2/d} & \mbox{for}\ d\geq 5, \end{cases}\end{equation}
		and $$t_{d,\alpha}:=\begin{cases} (4\alpha)^{-1}(4+((2\alpha+2d\alpha-d)\vee 0)) & \mbox{for}\ d<4,\\  (\alpha^{-1}\vee 7/2)-1 & \mbox{for}\ d=4,\\ 2(1+d^{-1}) & \mbox{for}\ d>4.\end{cases}$$ 
		The same bound holds for $|W_2^2(\tmu,\tnu)-W_2^2(\mu,\nu)|$ without assuming $\Tmn(\cdot)$ is Lipschitz. 
		
		\textbf{(Semi-discrete):} Set $\tmu=\mu$, $\tnu=\hnu$ or $\tmu=\hmu$, $\tnu=\nu$. Then the left hand side of~\eqref{eq:nsym} is $O_p(r_d^{(n,n)}\times (\log{(1+n)})^{t_{d,\alpha}})$ or $O_p(r_d^{(m,m)}\times (\log{(1+m)})^{t_{d,\alpha}})$ respectively.
	\end{theorem}
	
	A stronger result can be proved if both $\mu$ and $\nu$ are compactly supported.
	
	\begin{corollary}\label{cor:fsam}
		Consider the setting from~\cref{thm:nsmooth} and assume further that $\mu$ is compactly supported. Then, with $r_d^{(m,n)}$ defined as in~\eqref{eq:nsmrate}, we have:
		$$\E\left[\sup\limits_{\gamma\in\tgm} \int \lVert \ttg(\mx)-\Tmn(\mx)\rVert^2\,d\tmu(x)\right]\leq  C r^{(m,n)}_d,$$
		for some constant $C>0$, in both the discrete-discrete and semi-discrete settings from~\cref{thm:nsmooth}.
	\end{corollary}
	A brief description of the proof technique of~\cref{thm:nsmooth}  using~\cref{thm:newubd} is provided in~\cref{rem:pfnstech} below, and the actual proof is presented in~\cref{sec:mainrespf}.
	\begin{remark}[Proof technique]\label{rem:pfnstech} The proof of~\cref{thm:nsmooth} proceeds using the strategy outlined in~\cref{rem:pftech}. We first show that $\F_s$ (see~\cref{rem:pftech}) can be chosen as a certain class of convex functions which are in $L^2(\nu)$. We then use Dudley's entropy integral type bounds which in turn requires the bracketing entropy~\cite[Page 270]{vdV98} of $\F_s$, recently proved in~\cite[Equation 26]{kur2020convex}. This strategy is slightly different from that used in the proof of~\cite[Theorem 2]{chizat2020faster}, where the authors assume that $\mu$ is compactly supported whereas we only assume the finiteness of $\E \exp(t\lVert X_1\rVert^{\alpha})$ for some $t>0$, $\alpha>0$. The compactness assumption  on $\mu$ allows one to further restrict $\F_s$ to the class of Lipschitz functions. This additional restriction does not seem to be immediate without the compactness assumption.
	\end{remark}			
	As discussed in~\cref{sec:contrib}, the exponents obtained in~\cref{thm:nsmooth} are minimax optimal under bare minimal smoothness assumptions (see~\cite[Theorem 6]{hutter2021minimax}). To the best of our knowledge, rates for the discrete-discrete case for $m\neq n$ and those for the semi-discrete case were not known previously in the literature. Our rates are also strictly better than those (for different estimators, based on space tessellations) obtained in~\cite{berman2020convergence,li2020quantitative} and require less stringent assumptions than those in~\cite{chizat2020faster}. In the next section, we show how smoothness assumptions can be leveraged to mitigate the curse of dimensionality in~\cref{thm:nsmooth}.
	\subsection{Smooth plug-in estimator: mitigating the curse of dimensionality}\label{sec:smooth}
	
	In this section, we focus on two types of plug-in estimators for the densities associated with the probability measures $\mu$ and $\nu$: (a) wavelet based estimators (see~\cite{weed2019estimation,Tsybakov2009,Donoho1996,Kerk1992,Walter1992}) in~\cref{sec:wavelet}, and (b) kernel based estimators (see~\cite{Gine2008,Gine2016,Parzen1962,Silverman1978,Nadaraja1965}) in~\cref{sec:kernel}. In both these cases, we will show,  using~\cref{thm:newubd}, that the corresponding estimators of $\Tmn(\cdot)$ achieve (near) dimension-free rates under sufficient smoothness assumptions. 
	\subsubsection{Wavelet based estimators}\label{sec:wavelet}
	
	We begin this subsection by defining the Besov class of functions which will play a pivotal role in the sequel. 
	\begin{defn}[Besov class of functions]\label{def:besov}
		We describe Besov classes following the notation from~\cite[Section 2.1.1]{weed2019estimation}. Suppose $s>0$ and let $n>s$ be a positive integer. Given $\Omega\subseteq \R^d$, $h\in\R^d$ and $f(\cdot):\R^d\to\R^d$, set
		$$\boldsymbol{\Delta}_h^1 f(x):=f(x+h)-f(x),$$
		$$\boldsymbol{\Delta}_h^k f(x):= \boldsymbol{\Delta}_h^1\left(\boldsymbol{\Delta}_h^{k-1} f\right)(x), \quad \forall\ 2\leq k\leq n,$$
		where these functions are defined on $\Omega_{h,n}:=\{x\in\Omega:x+nh\in\Omega\}$. For $t>0$, we then define
		$$\omega_n(f,t):=\sup_{\lVert h\rVert\leq t} \lVert \boldsymbol{\Delta}_h^n f\rVert_{L^2(\Omega_{h,n})}.$$
		Finally, we define the space $\B^s(\Omega)$ to be the set of functions for which the quantity
		$$\lVert f\rVert_{\B^s(\Omega)}:=\lVert f\rVert_{L^2(\Omega)}+\sum_{j\geq 0} 2^{sj}\omega_n(f,2^{-j})$$
		is finite. The above expression can also be used to define Besov spaces (and norms) for $s<0$; see~\cite[Theorem 3.8.1]{Cohen2003}.
	\end{defn}
	
	In this subsection, we assume that $\mu$ and $\nu$ admit Besov smooth densities $\fm(\cdot)$ and $\fn(\cdot)$ (see~\cite{Cohen2003}~and~\cref{def:besov} above for details). Given $\Omega\subseteq \R^d$ and $s>0$, let $\B^s(\Omega)$ denote the set of Besov smooth functions on $\Omega$ of order $s$.
	\begin{Assumption}[Regularity of the densities]\label{as:besovden}
		Suppose that: \begin{enumerate} \item $\fm$ and $\fn$ are supported on compact and convex subsets of $\R^d$, say $\X$ and $\Y$ respectively. \item There exists $s,M> 0$ such that $\lVert\fm\rVert_{\B^{s}(\X)}\leq M$, $\lVert\fn\rVert_{\B^{s}(\Y)}\leq M$ and $\fm(x),\fn(y)\geq M^{-1}$ for all $x\in\X$, $y\in\Y$. \end{enumerate}
	\end{Assumption}
	We now present our wavelet based estimators for $\fm(\cdot)$ and $\fn(\cdot)$. Towards this direction, we begin with \emph{sets} of functions in $L^2(\X)$ (set of square integrable functions on $\X$), $\bph$ and $\{\bps_j\}_{j\geq 0}$, which form an orthonormal basis of $L^2(\X)$ and satisfy the standard regularity assumptions for a wavelet basis (see~\cite{Meyer1990,Hardle1998},~\cite[Appendix E]{weed2019estimation}). We defer a formal discussion on these assumptions to~\cref{def:wavelet} in the Appendix so as not to impede the flow of the paper. For the moment, it is worth noting that such sets of functions (e.g., Haar wavelets, Daubechies wavelets) are readily available in standard statistical softwares, see e.g., the \texttt{R} package \texttt{wavelets}.
	
	Next, fix $J_m\in\mathbb{N}$ (a truncation parameter to be chosen later depending on the sample size $m$). Consider the following:
	\begin{equation}\label{eq:estwave1}
	\hfm(x):=\sum_{\phi\in\bph}  a_{\phi}\phi(x)+\sum_{j=0}^{J_m} \sum_{\psi\in\bps_j} b_{\psi}\psi(x),
	\end{equation}
	where
	$$a_{\phi}:=\frac{1}{m}\sum_{i=1}^m \phi(X_i), \qquad b_{\psi}:=\frac{1}{m}\sum_{i=1}^m \psi(X_i).$$
	Unfortunately $\hfm(\cdot)$ as defined in~\eqref{eq:estwave1} may not be a probability density and consequently cannot be used to obtain plug-in estimators for $\ttg(\cdot)$. We therefore take the same route as in~\cite[Section 4.1]{weed2019estimation} to define the following estimator for $\fm(\cdot)$:
	\begin{equation}\label{eq:estwave2}
	\tfm:=\min_{g\in \mathcal{D}(\X)}\lVert g-\hfm\rVert_{\B^{-1}(\X)}, 
	\end{equation}
	where $\mathcal{D}(\X)$ is the space of probability density functions on $\X$ and $\B^{-1}(\X)$ is the Besov norm on $\X$ of order $-1$ as stated in~\cref{def:besov}. We can define $\tfn(\cdot)$ similarly. Computing both $\tfm(\cdot)$ and $\hfm(\cdot)$ (as it involves infinite sums) is challenging and we would refer the interested reader to~\cite[Section 6]{weed2019estimation} and the references therein, for details. Further  discussion of this aspect is beyond the scope of this paper.
	
	We are now in a position to present the main theorem of this subsection.
	\begin{theorem}\label{thm:smoothwav} Suppose that $\Tmn(\cdot)$ is $L$-Lipschitz, and $\tmu$ and $\tnu$ are the probability measures corresponding to the probability densities $\tfm(\cdot)$ and $\tfn(\cdot)$ with $m^{\frac{1}{d+2s}}\leq 2^{J_m}\leq m^{\frac{1}{d}}$ and $n^{\frac{1}{d+2s}}\leq 2^{J_n}\leq n^{\frac{1}{d}}$, then the following holds for some constant $C>0$:
		\begin{equation}\label{eq:smoothwav1}\E\left[\sup\limits_{\gamma\in\tgm} \int \lVert \ttg(\mx)-\Tmn(\mx)\rVert^2\,d\tmu(x)\right]\leq C \tilde{r}^{(m,n)}_{d,s},\end{equation}  \begin{equation}\label{eq:smoothwav2}\mbox{where}\quad \tilde{r}^{(m,n)}_{d,s}:=\begin{cases} m^{-1/2}\log{(1+m)}+n^{-1/2}\log{(1+n)} & \mbox{for}\ d=2,\\ m^{-\frac{1+s}{d+2s}}+n^{-\frac{1+s}{d+2s}} & \mbox{for}\ d\geq 3, \end{cases}\end{equation}	
		The same bound also holds for $\E|W_2^2(\tmu,\tnu)-W_2^2(\mu,\nu)|$.		
	\end{theorem}
	Note that $\frac{1+s}{d+2s}\to \frac{1}{2}$ as $s\to\infty$. Therefore~\cref{thm:smoothwav} shows that, when $m=n$, the rate of convergence for the wavelet based estimator is ``close" to $n^{-1/2}$ provided $s$ is large enough for each fixed $d$. This shows that $\Tmn(\cdot)$ obtained using the wavelet estimators for $\fm(\cdot)$ and $\fn(\cdot)$ mitigates the curse of dimensionality, contrast this with the  estimator in~\cref{thm:barate}. To avoid repetition, we defer further discussions on the rates observed in~\cref{thm:smoothwav} to~\cref{rem:curdim} where a holistic comparison is drawn with two other ``smooth'' plug-in estimators. 
	
	\subsubsection{Kernel based estimators}\label{sec:kernel}
	
	We first introduce the Sobolev class of functions which we will exploit in this subsection to construct estimators that achieve rates of convergence which mitigate the curse of dimensionality under sufficient smoothness.
	\begin{defn}[Uniform Sobolev class of functions]\label{def:holder}
		Let $\Omega\subseteq \R^d$ and $f(\cdot)$ be uniformly continuous on $\Omega$ and admits uniformly continuous derivatives up to order $s$ on $\Omega$ for some $s\in\mathbb{N}$. For any $\mm:=(m_1,\ldots ,m_d)\in\mathbb{N}^d$, let 
		$$\partial^{\mm} f:=\frac{\partial}{\partial_{x_1}^{m_1}}\ldots \frac{\partial}{\partial_{x_d}^{m_d}} f,\quad |\mm|:=\sum_{i=1}^d m_i.$$
		For any $k\leq s$, we further define,
		$$\lVert f\rVert_{C^k(\Omega)}:=\sum_{|\mm|\leq k} \lVert \partial^{\mm} f\rVert_{L^{\infty}(\Omega)}.$$
		
		The space $C^{s}(\Omega)$ is defined as the set of functions $f(\cdot)$ for which $\lVert f\rVert_{C^k(\Omega)}<\infty$ for all $k\leq s$.
	\end{defn}
	
	For this subsection, assume that $\mu$ and $\nu$ admit Sobolev smooth densities $\fm(\cdot)$ and $\fn$ in the uniform norm (see~\cref{def:holder} above). Given  $\Omega\subseteq \R^d$ and $s\in\mathbb{N}$, let $C^s(\Omega)$ denote the set of Sobolev smooth functions on $\Omega$ of order $s$. 
	
	\begin{Assumption}[Regularity of the densities]\label{as:smdensities}
		Suppose that \begin{enumerate} \item $\fm$ and $\fn$ are supported on compact and convex subsets of $\R^d$, say $\X$ and $\Y$ respectively. \item There exists $s,M> 0$ such that $\fm(\cdot)\in C^{s}(\X;M)$ and $\fn(\cdot)\in C^{s}(\Y;M)$ where $C^{s}(\X;M)$ is the space of real valued functions supported on $\X$ such that for all $f(\cdot)\in C^s(\X;M)$, we have $M^{-1}\leq f(x)\leq M$ for all $x\in\X$ and $\lVert f\rVert_{C^s(\X)}\leq M$. Here $\lVert\cdot\rVert_{C^s(\X)}$ is the standard uniform Sobolev norm as defined in~\cref{def:holder}. The space $C^s(\Y;M)$ is defined analogously. \end{enumerate}
	\end{Assumption}
	
	We now define our estimators for $\fm(\cdot)$ and $\fn(\cdot)$ using the standard kernel density estimation technique (see~\cite[Section 1.2]{Tsybakov2009}). Set
	\begin{equation}\label{eq:den}
	\hfm(x):=\frac{1}{mh_m^d}\sum_{i=1}^m K_d\left(\frac{X_i-x}{h_m}\right),
	\end{equation}
	for some bandwidth parameter $h_m>0$ and $d$-variate kernel $K_d(\cdot)$. We assume that $K_d(\cdot)$ is the $d$-fold product of  univariate kernels, i.e., there exists a kernel $K(\cdot)$ such that for $u=(u_1,\ldots ,u_d)\in\R^d$, $K_d(u)=\prod_{i=1}^d K(u_i)$. We define $\hfn(\cdot)$ similarly with the same univariate kernel and bandwidth.
	\begin{Assumption}[Regularity of the kernel]\label{as:kernel}
		Assume that $K(\cdot)$ is a symmetric, bounded, $s+1$ times differentiable kernel on $\R^d$ with all $s+1$ derivatives bounded and integrable. Further, suppose that $K(\cdot)$ is of order $2s+2$, i.e.,
		$$\int u^j K(u)\,du=\mathbbm{1}(j=0),\qquad\mbox{for }\  j=\{0,1,2,\ldots ,2s+1\},\qquad \mbox{and} \ \int |u|^{2s+2}|K(u)|\,du<\infty.$$				
	\end{Assumption}
	The above assumptions on $K(\cdot)$ are standard for estimating smooth densities and their derivatives of different orders in the kernel density estimation literature; see e.g.~\cite{Hansen08,arias2016estimation,Gine2008,Tsybakov2009,Gine2016}. There are several natural ways to construct kernels satisfying~\cref{as:kernel}, see~\cite[Section 1.2.2]{Tsybakov2009}; an example is also provided in~\cref{ex:ordker} below.
	
	\begin{example}[Example of a kernel satisfying~\cref{as:kernel}]\label{ex:ordker}
		Let $\psi_m(\cdot)$ be the $m$-th Hermite polynomial on $\R$ (see~\cite{laplace1820theorie}). Then the kernel function defined as
		$$K(u):=\sum_{m=0}^{2s+2} \psi_{m}(0)\psi_m(u)\exp(-u^2/2)$$
		satisfies~\cref{as:kernel}.
	\end{example}
	
	It is evident from~\cref{as:kernel} that $K(\cdot)$ may take some negative values, in which case, $\hfm(\cdot)$ (respectively $\hfn(\cdot)$) may not be a probability density. Consequently the barycentric projection (see~\cref{def:bcenterproj}) between $\hfm(\cdot)$ and $\hfn(\cdot)$ is not well-defined. We get around this by resorting to the same (approximate) projection technique we used in~\eqref{eq:estwave2} for the wavelet based estimators. In this case however, instead of (approximately) projecting using an appropriate Besov norm as in~\eqref{eq:estwave2}, we use a integral probability metric (see~\cref{def:ipm}; also see~\cite{Sriperumbudur2012,Alfred1997,rachev1985monge} for examples, computational procedures and applications of such metrics). The corresponding measure is defined below:
\begin{defn}[Integral probability metric]\label{def:ipm}
	Given a class $\F$ of bounded functions on $\R^d$ and two probability densities $g_1(\cdot)$ and $g_2(\cdot)$ on $\R^d$, the integral probability metric/distance between $g_1(\cdot)$ and $g_2(\cdot)$ with respect to $\F$ is defined as
	$$d_{\mathrm{IP}}(g_1,g_2;\F):=\sup_{\psi(\cdot)\in\F} \bigg|\int \psi(x)(g_1(x)-g_2(x))\,dx\bigg|.$$
	Sufficient conditions on $\F$ for $d_{\mathrm{IP}}(\cdot,\cdot;\F)$ to be a metric on the space of probability measures (not on the space of probability densities as they can be altered on set of Lebesgue measure $0$ without altering the underlying probability measures) on $\R^d$ have been discussed in~\cite{Alfred1997}. Observe that the measure $d_{\mathrm{IP}}(g_1,g_2;\F)$ is well defined even when $g_1(\cdot)$ and $g_2(\cdot)$ are not probability densities.
	
	In~\cref{thm:smooth} below, we use $\F=C^{s+2}(\X,M')$. Note that any function in $C^{s+2}(\X,M')$ can be extended to a function in $C^{s+2}(\R^d;M')$ (see~\cite[Theorem 23]{hutter2021minimax}~and~\cite[Theorem 1.105]{Triebel2006}). The fact that this choice of $\F$ results in a metric follows from the argument in~\cite[Page 8]{Alfred1997}. 
\end{defn}
 
 We are now in a position to describe the projection estimators for $f_{\mu}(\cdot)$ and $f_{\nu}(\cdot)$, and the rates achieved by the corresponding plug-in estimator. 
	\begin{theorem}\label{thm:smooth} Assume that $\Tmn(\cdot)$ is $L$-Lipschitz and $\fm$, $\fn$ are Lebesgue densities satisfying~\cref{as:smdensities}. Also suppose that $K(\cdot)$ satisfies~\cref{as:kernel}. Define $h_m:=m^{-\frac{1}{d+2s}}\log{m}$ , $h_n:=n^{-\frac{1}{d+2s}}\log{n}$ and $T:=\int |K_d(u)|\,du+1$. Fix any $M'>0$. Consider any probability density $\tfm^{M'}(\cdot)\in C^s(\X;TM)$ (where $M$ is defined as in~\cref{as:smdensities}) which satisfies
		\begin{equation}\label{eq:kerproj}
		d_{\mathrm{IP}}\left(\tfm^{M'},\hfm;C^{s+2}(\X;M')\right)\leq \inf_{\substack{f(\cdot)\in C^s(\X;TM)\\ f\geq 0,\ \int f=1}} d_{\mathrm{IP}}\left(\hfm,f;C^{s+2}(\X;M')\right)+r_{d,s}^{(m,n)}
		\end{equation}
		where $r_{d,s}^{(m,n)}$ is defined as in~\eqref{eq:kerproj1} and $d_{\mathrm{IP}}(\cdot,\cdot;C^{s+2}(\X;M'))$ is the integral probability metric defined in~\cref{def:ipm}. We define  $\tfn^{M'}(\cdot)$ analogously as in~\eqref{eq:kerproj} with $\X$, $\hfm(\cdot)$ replaced by $\Y$, $\hfn(\cdot)$. Then the following conclusions hold.
		
		\begin{enumerate}			
			\item There exists $M'>0$ (depending on $M$) such that, if, $\tmu$ and $\tnu$ are the probability measures corresponding to the probability densities $\tfm^{M'}(\cdot)$ and $\tfn^{M'}(\cdot)$, then the following holds for some constant $C>0$:
			$$\E\left[\sup\limits_{\gamma\in\tgm} \int \lVert \ttg(\mx)-\Tmn(\mx)\rVert^2\,d\tmu(x)\right]\leq C r_{d,s}^{(m,n)},$$  \begin{equation}\label{eq:kerproj1} \mbox{where}\quad r_{d,s}^{(m,n)} := \begin{cases} m^{-1/2}+n^{-1/2} & \mbox{for}\ d<2(s+2),\\ m^{-1/2}\left(\log{(1+m)}\right)^d+n^{-1/2}\left(\log{(1+n)}\right)^d & \mbox{for}\ d=2(s+2),\\ m^{-\frac{s+2}{d}}+n^{-\frac{s+2}{d}} & \mbox{for}\ d\geq 2(s+2). \end{cases}\end{equation}
			The same bound also holds for $\E|W_2^2(\tmu,\tnu)-W_2^2(\mu,\nu)|$.
			
			\item $\hfm(\cdot)$ satisfies
			\begin{equation}\label{eq:kerproj2}
			\lim\limits_{n\to\infty} \max\left\{\P\left(\lVert \hfm\rVert_{C^s(\tilde{\X})}\geq TM\right),\P\left(\sup_{x\in\tilde{\X}} |\hfm(x)-\fm(x)|\geq \eps\right)\right\}=0
			\end{equation}
			for any $\eps>0$, where $\tilde{\X}$ is any compact subset of $\X^o$. The same conclusion holds for $\hfn(\cdot)$ with $\X$ replaced by $\Y$.
		\end{enumerate}
	\end{theorem}
	In~\cref{thm:smooth}, $\tfm^{M'}(\cdot)$ can be viewed as an approximate minimizer of $d_{\mathrm{IP}}(\hfm,\cdot;C^{s+2}(\X,M'))$ over an appropriate class of Sobolev smooth probability densities. This is carried out because $\hfm(\cdot)$ by itself may not be a probability density.
	
	Further note that $\tmu,\tnu$ as specified in~\cref{thm:smooth} are both smooth, and consequently $\tgm$ is a singleton and the supremum in~\cref{thm:smooth} can be dropped. A brief description of the proof technique for~\cref{thm:smooth} is presented in~\cref{rem:pfstech} below  and the actual proof is given in~\cref{sec:mainrespf}. 
	\begin{remark}[Proof technique]\label{rem:pfstech} The proof of~\cref{thm:smooth} proceeds along the same lines as~\cref{rem:pfnstech}. We first show that $\F_s$ (see~\cref{rem:pftech}) can be chosen as a certain subset of $C^{s+2}(\Y^{\circ})$. We then use Dudley's entropy integral type bounds which in turn requires the bracketing entropy~\cite[Page 270]{vdV98} of the class of compactly supported Sobolev smooth functions which can be found in~\cite[Corollary 2.7.2]{vaart2013}.
	\end{remark}
	We now explain the implications of both the parts of~\cref{thm:smooth} in the following two remarks.
	
	\begin{remark}[Mitigating the curse of dimensionality]\label{rem:curdim}
		~\cref{thm:smooth} shows that, under enough smoothness, i.e., when $2(s+2)>d$, both the upper bounds for~\eqref{eq:OTrate}~and~\eqref{eq:wassrate} are $O_p(n^{-1/2})$. This shows that, for large dimensions, provided $\mu$ and $\nu$ admit smooth enough densities, it is possible to construct plug-in estimators that mitigate the curse of dimensionality. Note that a similar estimator was analyzed in~\cite[Proposition 1]{gunsilius2018convergence} when $m=n$. However, the rates obtained in~\cref{thm:smooth} are \emph{strictly better} than those in~\cite[Proposition 1]{gunsilius2018convergence}. For $m=n$, when $d<2(s+2)$,~\cite{gunsilius2018convergence} obtained a rate of $n^{-\frac{s+2}{2(s+2)+d}}$ which is \emph{worse} than $n^{-1/2}$ obtained in~\cref{thm:smooth}. For the other regimes,~\cite{gunsilius2018convergence} obtains rates (up to log factors) of $n^{-1/4}$ and $n^{-\frac{1}{(s+2)(d+2(s+2))}}$ which are both worse than the respective rates of $n^{-1/2}$ and $n^{-\frac{s+2}{d}}$ in~\cref{thm:smooth}. In fact, the rates obtained in~\cref{thm:smoothwav} are also strictly better than the rates obtained in~\cite[Proposition 1]{gunsilius2018convergence} (for every fixed $d$) described above, but they are strictly worse than the rates obtained in~\cref{thm:smooth}. When the degree of smoothness $s$ is  large, both Theorems~\ref{thm:smoothwav}~and~\ref{thm:smooth} lead to rates of (approximately) $n^{-1/2}$. However when $s$ is small, the rate in~\cref{thm:smooth} is much faster than that in~\cref{thm:smoothwav}, e.g., if $s$ is close to $0$, the rate in~\cref{thm:smoothwav} is approximately $n^{-\frac{1}{d}}$ whereas that in~\cref{thm:smooth} is the faster rate of $n^{-\frac{2}{d}}$. It must be noted however that the smoothness assumptions are different in Theorems~\ref{thm:smoothwav}~and~\ref{thm:smooth}.
	\end{remark}
	
	\begin{remark}[Computational aspects of \cref{thm:smooth}]\label{rem:compasp}
		In~\cref{thm:smooth}, we have shown that the plug-in estimator for $\Tmn(\cdot)$ using $\tfm^{M'}(\cdot)$ and $\tfn^{M'}(\cdot)$ achieve rates that mitigate the curse of dimensionality under sufficient smoothness. However, as is evident, $\tfm^{M'}(\cdot)$ is hard to compute whereas $\hfm(\cdot)$ is computable easily in linear time. Note that if $\hfm(\cdot)$ itself were a probability density in $C^s(\X;TM)$, then we would have $\hfm=\tfm^{M'}$. While~\cref{thm:smooth} does not establish that, it does come close in part 2, from which we can easily derive the following:
		$$\lim_{n\to\infty}\P(\hfm(\cdot)\notin C^s(\tilde{\X};TM))=0.$$
		The above shows that $\hfm(\cdot)$ is indeed bounded below by $(TM)^{-1}$ on $\tilde{\X}$ (any compact subset of the interior of $\X$), and additionally belongs to $C^s(\tilde{\X};TM)$ with probability converging to $1$. This leads us to conjecture that the natural density version of $\hfm(\cdot)$, i.e.,
		$$\frac{\max\{\hfm(\cdot),0\}}{\int \max\{\hfm(x),0\}\,dx}$$
		should serve as a good proxy for $\tfm^{M'}(\cdot)$ and lead to rates of convergence that mitigate the curse of dimensionality. From a computational perspective, the density specified above is easy to simulate from using an accept-reject algorithm without computing the integral in the denominator (see~\cite[Algorithm 4.3]{oudjane2005sup}). However, our current proof technique does not provide rates of convergence for the above density estimator based on $\hfm(\cdot)$.
	\end{remark}
	Another important  implication of~\cref{thm:smooth} is the bound obtained on $|W_2(\tmu,\tnu)-W_2(\mu,\nu)|$ when $\mu\neq\nu$. We first present the result and then describe the implication.
	
	\begin{proposition}\label{prop:wassrate}
		Consider the setting in~\cref{thm:smooth}. Then, provided $\mu\neq \nu$, the following holds:
		$$|W_2(\tmu,\tnu)-W_2(\mu,\nu)|=O_p(r_{d,s}^{(m,n)}).$$
	\end{proposition}
	~\cref{prop:wassrate} (see~\cref{sec:mainrespf} for a proof) shows an interesting  distinction between the $\mu\neq\nu$ case and the $\mu=\nu$ case. For $\mu=\nu$, the best possible exponent is $n^{-\frac{1+s}{2s+d}}$ for $d\geq 3$ (see~\cite[Theorem 3]{weed2019estimation} where the result was established under more general Besov smoothness assumptions). On the contrary, when $\mu\neq\nu$,~\cref{prop:wassrate} establishes a rate of $n^{-\frac{s+2}{d}}$ for the Wasserstein distance which is \emph{strictly better} than the minimax achievable rate mentioned above when $\mu=\nu$. This observation complements~\cite[Corollary 1]{chizat2020faster} where the authors make a similar observation for the special case of $s=0$.
	\subsection{Discretized plug-in estimator under smoothness assumptions}\label{sec:dismooth}
	In~\cref{sec:nsmooth}, we discussed how smoothness can be incorporated into the plug-in procedure to get faster rates of convergence. Such plug-in estimators are popular in the computational OT literature (see~\cite{benamou2000computational,benamou2014numerical,chartrand2009gradient,cuturi2013sinkhorn}). However, even after $\tfm(\cdot)$, $\tfn(\cdot)$ are calculated, $\ttg$ as in~\cref{thm:smooth} cannot be computed explicitly from data if $\tfm(\cdot)$ and $\tfn(\cdot)$ are continuous densities. This is in contrast to $\ttg$ from~\cref{thm:nsmooth} in the \emph{discrete-discrete} case which is explicitly computable using a standard linear program, but achieves worse rates of convergence. This is not unexpected. Thanks to the \emph{no free lunch} principle, better statistical accuracy is naturally accompanied by heavier computational challenges. Therefore, our goal here is to construct estimators, under smoothness assumptions as in~\cref{sec:smooth}, which are computable in polynomial time (with complexity increasing with smoothness) provided $\tfm(\cdot)$ and $\tfn(\cdot)$ can be sampled from, and also attain rates that mitigate the curse of dimensionality.
	
	\emph{Construction}: We will illustrate the discretized estimator using the kernel based estimator from~\cref{sec:kernel}. Similar results also hold for the wavelet based estimator from~\cref{sec:wavelet}. Recall the kernel density estimators $\tfm(\cdot)$ and $\tfn(\cdot)$ (see~\eqref{eq:kerproj}). Sample $M\geq 1$ random points from both $\tfm(\cdot)$ and $\tfn(\cdot)$. Let $\tM$ and $\tN$ denote the standard empirical measures on the $M$ points sampled from $\tfm(\cdot)$ and $\tfn(\cdot)$ respectively. Finally construct $\Tes\equiv \ttg$ as in~\cref{def:bcenterproj} with $\tmu=\tM$ and $\tnu=\tN$. It should be pointed out that a similar construction was also used in~\cite[Section 6]{weed2019estimation} for estimating probability densities under the Wasserstein loss. Based on this construction, the main result of this section is as follows:
	
	\begin{theorem}\label{thm:dismooth}
		Consider the setting in~\cref{thm:smooth} and the same construction of $\ttg$ as above. For simplicity, let's also assume $m=n$. Accordingly set $M=n^{\frac{s+2}{2}}$. Then $\tgm$ is a singleton and consequently the following conclusion holds for some constant $C>0$:
		$$\E\left[\int \lVert \Tes(\mx)-\Tmn(\mx)\rVert^2\,d\tmu(x)\right]\leq C r_{d,s}^{(n,n)}.$$
		The same rates also hold for $\E|W_2^2(\tmu,\tnu)-W_2^2(\mu,\nu)|$.
	\end{theorem}
	The proof of~\cref{thm:dismooth} is given in~\cref{sec:mainrespf}. Once the empirical measures $\tM$ and $\tN$ have been obtained, an explicit computation of $\Tes$ as described above requires  $O(M^3)=O(n^{\frac{3(s+2)}{2}})$ steps using the \emph{Hungarian algorithm}, see~\cite{jonker1987shortest}. This highlights the statistical versus computational trade off, i.e., in order to mitigate the curse of dimensionality in convergence rates by exploiting smoothness, the computational complexity gets progressively worse by polynomial factors in $n$. It should be mentioned that (approximate) algorithms faster than the Hungarian algorithm stated above, can be found in~\cite{gabow1989faster,agarwal2014approximation,cuturi2013sinkhorn} to name a few. Due to space constraints, we avoid a detailed discussion on this. 
	
	In the above construction, sampling from the smoothed kernel densities $\tfm(\cdot)$ and $\tfn(\cdot)$ is crucial. If we would simply draw $M$ bootstrap samples from the empirical distributions $\hmu$ and $\hnu$, the rates of convergence wouldn't improve from those observed in~\cref{thm:nsmooth} no matter how large $M$ is chosen.
	\section{Applications}\label{sec:App}
	In this section, we will apply our results to two popular problems, namely --- estimating the Wasserstein barycenter between two probability distributions (see~\cite{agueh2011barycenters,cuturi2014fast,bigot2018characterization,carlier2015numerical}) in~\cref{sec:bar}, and obtaining detection thresholds in some recent optimal transport based independence testing procedures (see~\cite{deb2021multivariate,deb2020measuring,ghosal2019multivariate,shi2020distribution,Shi2020}) in~\cref{sec:nonptest}.
	\subsection{Wasserstein barycenter estimation}\label{sec:bar}
	Let $\mu,\nu\in\Ptac$. The Wasserstein barycenter between $\mu$ and $\nu$ is then given by:
	\begin{equation}\label{eq:barycenter}
	\Bmn:=\min\limits_{\rho\in\Pac} \left(\frac{1}{2}W_2^2(\mu,\rho)+\frac{1}{2}W_2^2(\rho,\nu)\right).
	\end{equation}
	In fact, by~\cref{prop:bmopt}, there exists an optimal transport map $T_0$ from $\mu$ to $\nu$ and by~\cite{agueh2011barycenters,bigot2018characterization,boissard2015distribution}, an alternative characterization of $\Bmn$ is as follows:
	\begin{equation}\label{eq:altrep}
	\Bmn=\left(\frac{1}{2}\mbox{Id}+\frac{1}{2}\Tmn\right)\#\mu,\qquad \mbox{where}\qquad \mbox{Id}(x)=x.
	\end{equation}
	Estimating $\Bmn$ as in~\eqref{eq:barycenter} has attracted significant attention over the past few years in economics~\cite{chiappori2010hedonic,carlier2010matching}, Bayesian learning~\cite{srivastava2018scalable,srivastava2015wasp}, dynamic formulations~\cite{claici2018stochastic,chewi2020gradient}, algorithmic fairness~\cite{gordaliza2019obtaining,chzhen2019}, etc. The most natural strategy employed in estimating $\Bmn$ is to use the empirical plug-in estimator, i.e., replacing $\mu,\nu$ in~\eqref{eq:barycenter} with $\hmu,\hnu$. This strategy has been used, approximated and analyzed extensively in e.g.,~\cite{cuturi2014fast,carlier2015numerical,le2017existence,boissard2015distribution}. Based on~\eqref{eq:altrep}, the natural plug-in estimator of $\Bmn$ would be:
	\begin{equation}\label{eq:estbar}
	\Bmnz^{\gamma}=\left(\frac{1}{2}\mbox{Id}+\frac{1}{2}\ttg\right)\#\tmu
	\end{equation}
	where $\ttg$ is the plug-in estimator of $T_0$ obtained by solving~\eqref{eq:conmean}, with $\mu$ and $\nu$ replaced by $\tmu$ and $\hnu$ respectively and $\gamma\in\tgm$. While the consistency of $\Bmnz^{\gamma}$ has been analyzed for $m=n$  in~\cite{le2017existence} and rates have been obtained for $d=1$ in~\cite{bigot2018upper}, the more general question of obtaining rates of convergence for $\Bmnz^{\gamma}$ for general dimensions $d\geq 1$ is yet unanswered. We address this question in the following result (see~\cref{sec:Appf} for a proof).
	\begin{theorem}\label{thm:barate}
		Suppose that the same assumptions from~\cref{thm:nsmooth} hold. Then, with $\Bmnz^{\gamma}$ as defined in~\eqref{eq:estbar} and $r_d^{(m,n)}$, $t_{d,\alpha}$ defined in~\cref{thm:nsmooth}, the following holds:
		$$\sup\limits_{\gamma\in\tgm}W_2^2(\Bmnz^{\gamma},\Bmn)=O_p\left(r^{(m,n)}_d\times (\log{(1+\max\{m,n\})})^{t_{d,\alpha}}\right).$$
	\end{theorem}
	\subsection{Nonparametric independence testing: Optimal transport based Hilbert-Schmidt independence criterion}\label{sec:nonptest}
	Let $(X_1,Y_1),\ldots ,(X_n,Y_n)\overset{i.i.d.}{\sim}\pi$, a probability measure on $\R^{d_1+d_2}$, with marginals $\mu\in\Paco$ and $\nu\in\Pact$. Our problem of interest is the following hypothesis testing problem, given as:
	\begin{equation}\label{eq:indprob}
	\Hy_0:\pi=\mu\otimes \nu \qquad \mbox{versus} \qquad \Hy_1:\pi\neq \mu\otimes\nu.
	\end{equation}
	This is the classical nonparametric independence testing problem which has received a lot of attention in the statistics and machine learning literature (see~\cite{SzekelyCorrDist07,GrettonKernelMeasInd05,Heller2013,Berrett19}, and~\cite{Drouet2001,Josse2016} for a review). In keeping with the overall theme of this paper, our focus here will be on a large class of OT based  independence testing procedures, introduced first in~\cite{deb2021multivariate} followed by recent developments in~\cite{shi2020distribution,Shi2020,deb2020measuring}. These tests bear resemblance to the Hilbert-Schmidt independence criterion (HSIC); see~\cite{GrettonKernelMeasInd05,gretton2005measuring,gretton2007kernel} and have attractive properties such as distribution-freeness (see~\cref{prop:testprop}), consistency without moment assumptions and robustness against heavy-tailed distributions and against  contamination~\cite{deb2021multivariate,shi2020distribution}. Below, we describe this class of tests, see~\eqref{eq:emprhsic}~and~\eqref{eq:testrhsic}. Our main theoretical contribution of this section will be to provide detection thresholds of these OT based tests.
	
	\emph{Construction}: Suppose $\upsilon_1$, $\upsilon_2$ be two compactly supported probability distributions on $\R^{d_1}$ and $\R^{d_2}$ respectively (e.g., $\upsilon_1\equiv \U[0,1]^{d_1}$, $\upsilon_2\equiv \U[0,1]^{d_2}$). Let $U_1,\ldots ,U_n\overset{i.i.d.}{\sim}\upsilon_1$, $V_1,\ldots ,V_n\overset{i.i.d.}{\sim}\upsilon_2$,  $\hu:=n^{-1}\sum_{i=1}^n \delta_{U_i}$ and $\hv:=n^{-1}\sum_{j=1}^n \delta_{V_j}$. Recall the definitions of $\hmn$ (with $m=n$) and $\hnu$ from~\eqref{eq:empdef}. Let $\hto$ ($\htt$) be obtained by solving~\eqref{eq:conmean}, with $\mu$ and $\nu$ replaced by $\hmn$ and $\hu$ ($\hnu$ and $\hv$) respectively. Consider two non negative definite, continuous, characteristic kernels (see~\cite{fukumizu2008characteristic,Sriperumbudur2010} for definitions) $K_1(\cdot,\cdot)$ and $K_2(\cdot,\cdot)$ on $(\mbox{supp}(\upsilon_1))^2$ and $(\mbox{supp}(\upsilon_2))^2$. Set $\hx_{ij}:=K_1(\hto(X_i),\hto(X_j))$ and $\hy_{ij}:=K_2(\htt(Y_i),\htt(Y_j))$. Our test statistic is as follows:
	\begin{equation}\label{eq:emprhsic}\hhs:=n^{-2}\sum_{i,j}\hx_{ij}\hy_{ij}+n^{-4}\sum_{i,j,r,s}\hx_{ij}\hy_{rs}-2n^{-3}\sum_{i,j,r} \hx_{ij}\hy_{ir}.\end{equation}
	\begin{proposition}[See~\cite{deb2021multivariate,Shi2020}]\label{prop:testprop}
		\begin{enumerate}
			\item \textbf{(Distribution-freeness)} When $X_1$ and $Y_1$ are independent, the distribution of $n\times \hhs$ is universal, i.e., it does not depend on $\mu$ and $\nu$ for every fixed $n$.
			\item \textbf{(Consistency against fixed alternatives)} Let $c_{n,\alpha}$ be the upper $(1-\alpha)$-th quantile from the universal distribution in part 1 above. Then $\hhs\overset{P}{\longrightarrow}\nhs(\pi|\mu\otimes\nu)$ where 
			\begin{align}\label{eq:popmeas}
			&\nhs(\pi| \mu\otimes\nu):=\E [K_1(T_1(X_1),T_1(X_2))K_2(T_2(Y_1),T_2(Y_2))]+\E [K_1(T_1(X_1),T_1(X_2))]\nonumber \\&\times \E [K_2(T_2(Y_1),T_2(Y_2))]-2\E [K_1(T_1(X_1),T_1(X_2))K_2(T_2(Y_1),T_2(Y_3))],
			\end{align}
			where $T_{1}(\cdot)$ (respectively $T_2(\cdot)$) is the optimal transport map from $\mu$ ($\nu$) to $\upsilon_1$  ($\upsilon_2$); see~\cref{def:otmpm}. Further $\nhs(\pi|\mu\otimes\nu)=0$ if and only if $\pi=\mu\otimes\nu$. Define the following test function:
			\begin{equation}\label{eq:testrhsic}
			\phi_{n,\alpha}:=\mathbbm{1}(n\times\hhs\geq c_{n,\alpha}).
			\end{equation}
			Then $\E[\phi_{n,\alpha}]\to 1$ as $n\to\infty$ under $\Hy_1$, i.e., when $\pi\neq \mu\otimes\nu$.
		\end{enumerate}
	\end{proposition}
	~\cref{prop:testprop} shows that the test based on $\hhs$ (see~\eqref{eq:emprhsic}), i.e., $\phi_{n,\alpha}$ (see~\eqref{eq:testrhsic}), can be carried out without resorting to the permutation principle as is necessary for the usual HSIC based test (see~\cite{gretton2005measuring}). Further, when the sampling distribution is fixed,~\cref{prop:testprop} shows that $\hhs$ consistently estimates $\nhs(\pi|\mu\otimes\nu)$, a quantity which equals $0$ if and only if $\pi=\mu\otimes\nu$ (this yields the consistency of $\phi_{n,\alpha})$ against fixed alternatives. 
	
	While consistency against fixed alternatives is an attractive feature of $\phi_{n,\alpha}$, a more intricate question of statistical interest is to understand the local power of $\phi_{n,\alpha}$ under ``changing sequence of alternatives converging to the null" as $n\to\infty$. To study the local power of $\phi_{n,\alpha}$, we need to consider a triangular array setting, where the data distribution changes with $n$, i.e., $(X_1,Y_1),\ldots ,(X_n,Y_n)\overset{i.i.d.}{\sim}\pn$, a probability measure on $\R^{d_1+d_2}$, with marginals $\mn\in\Paco$ and $\nn\in\Pact$. As $\nhs(\cdot|\cdot)$ characterizes independence, a mathematical formulation of ``alternatives converging to null" would be to say $\nhs(\pn|\mn\otimes\nn)\to 0$ as $n\to\infty$. Similar questions have attracted a lot of attention in modern statistics, featuring measures (other than $\nhs(\cdot|\cdot)$) which characterize independence, see e.g.,~\cite{Berrett19,kim2020minimax,li2019optimality,auddy2021exact}. In the following result (see~\cref{sec:Appf} for a proof), we show that if $\nhs(\pn|\mn\otimes\nn)\to 0$ slowly enough with $n$, then $\phi_{n,\alpha}$ yields a consistent sequence of tests for problem~\eqref{eq:indprob}.
	
	\begin{theorem}\label{thm:dethresh}
		Consider problem~\eqref{eq:indprob} with $\pn$, $\mn$, $\nn$ (changing with $n$) and suppose $T_{1,n}(\cdot)$ and $T_{2,n}(\cdot)$ are both $L$-Lipschitz ($L$ is free of $n$). Also assume $K_1(\cdot)$, $K_2(\cdot)$ are Lipschitz, $\mn$, $\nn$ are supported on fixed compact sets (supports are free of $n$). Set $r_{d_1,d_2}^{(n,n)}:=r_{d_1}^{(n,n)}+r_{d_2}^{(n,n)}$ where $r_{d_1}^{(n,n)},r_{d_2}^{(n,n)}$ is defined via~\eqref{eq:nsmrate}. Then, 
		$$\E[\phi_{n,\alpha}]\to 1\qquad \mbox{if} \qquad (r_{d_1,d_2}^{(n,n)})^{-1/2}\times\nhs(\pn|\mn\otimes\nn)\to\infty,$$
	\end{theorem}

	\bibliographystyle{plain}
	\bibliography{OT}

\begin{thebibliography}{}

\bibitem[Agarwal and Sharathkumar, 2014]{agarwal2014approximation}
Agarwal, P.~K. and Sharathkumar, R. (2014).
\newblock Approximation algorithms for bipartite matching with metric and
  geometric costs.
\newblock In {\em Proceedings of the forty-sixth annual ACM symposium on Theory
  of computing}, pages 555--564.

\bibitem[Agueh and Carlier, 2011]{agueh2011barycenters}
Agueh, M. and Carlier, G. (2011).
\newblock Barycenters in the wasserstein space.
\newblock {\em SIAM Journal on Mathematical Analysis}, 43(2):904--924.

\bibitem[Ambrosio and Gigli, 2013]{Ambrosio2013}
Ambrosio, L. and Gigli, N. (2013).
\newblock A user's guide to optimal transport.
\newblock In {\em Modelling and optimisation of flows on networks}, volume 2062
  of {\em Lecture Notes in Math.}, pages 1--155. Springer, Heidelberg.

\bibitem[Ambrosio et~al., 2008]{Ambro08}
Ambrosio, L., Gigli, N., and Savar\'{e}, G. (2008).
\newblock {\em Gradient flows in metric spaces and in the space of probability
  measures}.
\newblock Lectures in Mathematics ETH Z\"{u}rich. Birkh\"{a}user Verlag, Basel,
  second edition.

\bibitem[Arias-Castro et~al., 2016]{arias2016estimation}
Arias-Castro, E., Mason, D., and Pelletier, B. (2016).
\newblock On the estimation of the gradient lines of a density and the
  consistency of the mean-shift algorithm.
\newblock {\em The Journal of Machine Learning Research}, 17(1):1487--1514.

\bibitem[Auddy et~al., 2021]{auddy2021exact}
Auddy, A., Deb, N., and Nandy, S. (2021).
\newblock Exact detection thresholds for chatterjee's correlation.
\newblock {\em arXiv preprint arXiv:2104.15140}.

\bibitem[Barthe and Bordenave, 2013]{Barthe2013}
Barthe, F. and Bordenave, C. (2013).
\newblock Combinatorial optimization over two random point sets.
\newblock In {\em S\'{e}minaire de {P}robabilit\'{e}s {XLV}}, volume 2078 of
  {\em Lecture Notes in Math.}, pages 483--535. Springer, Cham.

\bibitem[Benamou and Brenier, 2000]{benamou2000computational}
Benamou, J.-D. and Brenier, Y. (2000).
\newblock A computational fluid mechanics solution to the monge-kantorovich
  mass transfer problem.
\newblock {\em Numerische Mathematik}, 84(3):375--393.

\bibitem[Benamou et~al., 2014]{benamou2014numerical}
Benamou, J.-D., Froese, B.~D., and Oberman, A.~M. (2014).
\newblock Numerical solution of the optimal transportation problem using the
  monge--amp{\`e}re equation.
\newblock {\em Journal of Computational Physics}, 260:107--126.

\bibitem[Berman, 2020]{berman2020convergence}
Berman, R.~J. (2020).
\newblock Convergence rates for discretized monge--amp{\`e}re equations and
  quantitative stability of optimal transport.
\newblock {\em Foundations of Computational Mathematics}, pages 1--42.

\bibitem[Bernton et~al., 2017]{bernton2017inference}
Bernton, E., Jacob, P.~E., Gerber, M., and Robert, C.~P. (2017).
\newblock Inference in generative models using the {W}asserstein distance.
\newblock {\em arXiv preprint arXiv:1701.05146}, 1(8):9.

\bibitem[Berrett and Samworth, 2019]{Berrett19}
Berrett, T.~B. and Samworth, R.~J. (2019).
\newblock Nonparametric independence testing via mutual information.
\newblock {\em Biometrika}, 106(3):547--566.

\bibitem[Bigot et~al., 2017]{bigot2017geodesic}
Bigot, J., Gouet, R., Klein, T., L{\'o}pez, A., et~al. (2017).
\newblock Geodesic {PCA} in the {W}asserstein space by convex {PCA}.
\newblock In {\em Annales de l'Institut Henri Poincar{\'e}, Probabilit{\'e}s et
  Statistiques}, volume~53, pages 1--26. Institut Henri Poincar{\'e}.

\bibitem[Bigot et~al., 2018]{bigot2018upper}
Bigot, J., Gouet, R., Klein, T., Lopez, A., et~al. (2018).
\newblock Upper and lower risk bounds for estimating the wasserstein barycenter
  of random measures on the real line.
\newblock {\em Electronic journal of statistics}, 12(2):2253--2289.

\bibitem[Bigot and Klein, 2018]{bigot2018characterization}
Bigot, J. and Klein, T. (2018).
\newblock Characterization of barycenters in the wasserstein space by averaging
  optimal transport maps.
\newblock {\em ESAIM: Probability and Statistics}, 22:35--57.

\bibitem[Blanchet and Carlier, 2016]{blanchet2016optimal}
Blanchet, A. and Carlier, G. (2016).
\newblock Optimal transport and cournot-nash equilibria.
\newblock {\em Mathematics of Operations Research}, 41(1):125--145.

\bibitem[Boeckel et~al., 2018]{boeckel2018multivariate}
Boeckel, M., Spokoiny, V., and Suvorikova, A. (2018).
\newblock Multivariate brenier cumulative distribution functions and their
  application to non-parametric testing.
\newblock {\em arXiv preprint arXiv:1809.04090}.

\bibitem[Boissard et~al., 2015]{boissard2015distribution}
Boissard, E., Le~Gouic, T., Loubes, J.-M., et~al. (2015).
\newblock Distribution’s template estimate with wasserstein metrics.
\newblock {\em Bernoulli}, 21(2):740--759.

\bibitem[Bonneel et~al., 2011]{bonneel2011displacement}
Bonneel, N., Van De~Panne, M., Paris, S., and Heidrich, W. (2011).
\newblock Displacement interpolation using lagrangian mass transport.
\newblock In {\em Proceedings of the 2011 SIGGRAPH Asia Conference}, pages
  1--12.

\bibitem[Boucheron et~al., 2013]{Boucheron2013}
Boucheron, S., Lugosi, G., and Massart, P. (2013).
\newblock {\em Concentration inequalities}.
\newblock Oxford University Press, Oxford.
\newblock A nonasymptotic theory of independence, With a foreword by Michel
  Ledoux.

\bibitem[Brenier, 1991]{Brenier1991}
Brenier, Y. (1991).
\newblock Polar factorization and monotone rearrangement of vector-valued
  functions.
\newblock {\em Comm. Pure Appl. Math.}, 44(4):375--417.

\bibitem[Bronshtein, 1976]{bronshtein1976varepsilon}
Bronshtein, E.~M. (1976).
\newblock $\varepsilon$-entropy of convex sets and functions.
\newblock {\em Siberian Mathematical Journal}, 17(3):393--398.

\bibitem[Caffarelli, 1992a]{Caf1992}
Caffarelli, L.~A. (1992a).
\newblock Boundary regularity of maps with convex potentials.
\newblock {\em Comm. Pure Appl. Math.}, 45(9):1141--1151.

\bibitem[Caffarelli, 1992b]{Caf92}
Caffarelli, L.~A. (1992b).
\newblock The regularity of mappings with a convex potential.
\newblock {\em J. Amer. Math. Soc.}, 5(1):99--104.

\bibitem[Caffarelli, 1996]{Caf1996}
Caffarelli, L.~A. (1996).
\newblock Boundary regularity of maps with convex potentials. {II}.
\newblock {\em Ann. of Math. (2)}, 144(3):453--496.

\bibitem[Carlier and Ekeland, 2010]{carlier2010matching}
Carlier, G. and Ekeland, I. (2010).
\newblock Matching for teams.
\newblock {\em Economic theory}, 42(2):397--418.

\bibitem[Carlier et~al., 2015]{carlier2015numerical}
Carlier, G., Oberman, A., and Oudet, E. (2015).
\newblock Numerical methods for matching for teams and wasserstein barycenters.
\newblock {\em ESAIM: Mathematical Modelling and Numerical Analysis},
  49(6):1621--1642.

\bibitem[Chartrand et~al., 2009]{chartrand2009gradient}
Chartrand, R., Wohlberg, B., Vixie, K., and Bollt, E. (2009).
\newblock A gradient descent solution to the monge-kantorovich problem.
\newblock {\em Applied Mathematical Sciences}, 3(22):1071--1080.

\bibitem[Chernozhukov et~al., 2017]{Chernozhukhov2017}
Chernozhukov, V., Galichon, A., Hallin, M., and Henry, M. (2017).
\newblock Monge-{K}antorovich depth, quantiles, ranks and signs.
\newblock {\em Ann. Statist.}, 45(1):223--256.

\bibitem[Chewi et~al., 2020]{chewi2020gradient}
Chewi, S., Maunu, T., Rigollet, P., and Stromme, A.~J. (2020).
\newblock Gradient descent algorithms for bures-wasserstein barycenters.
\newblock In {\em Conference on Learning Theory}, pages 1276--1304. PMLR.

\bibitem[Chiappori et~al., 2010]{chiappori2010hedonic}
Chiappori, P.-A., McCann, R.~J., and Nesheim, L.~P. (2010).
\newblock Hedonic price equilibria, stable matching, and optimal transport:
  equivalence, topology, and uniqueness.
\newblock {\em Economic Theory}, 42(2):317--354.

\bibitem[Chizat et~al., 2018]{chizat2018scaling}
Chizat, L., Peyr{\'e}, G., Schmitzer, B., and Vialard, F.-X. (2018).
\newblock Scaling algorithms for unbalanced optimal transport problems.
\newblock {\em Mathematics of Computation}, 87(314):2563--2609.

\bibitem[Chizat et~al., 2020]{chizat2020faster}
Chizat, L., Roussillon, P., L{\'e}ger, F., Vialard, F.-X., and Peyr{\'e}, G.
  (2020).
\newblock Faster {W}asserstein distance estimation with the sinkhorn
  divergence.
\newblock {\em Advances in Neural Information Processing Systems}, 33.

\bibitem[Chzhen et~al., 2019]{chzhen2019}
Chzhen, E., Denis, C., Hebiri, M., Oneto, L., and Pontil, M. (2019).
\newblock {Leveraging Labeled and Unlabeled Data for Consistent Fair Binary
  Classification}.
\newblock In {\em {NeurIPS 2019 - 33th Annual Conference on Neural Information
  Processing Systems}}, Vancouver, Canada.

\bibitem[Claici et~al., 2018]{claici2018stochastic}
Claici, S., Chien, E., and Solomon, J. (2018).
\newblock Stochastic wasserstein barycenters.
\newblock In {\em International Conference on Machine Learning}, pages
  999--1008. PMLR.

\bibitem[Cohen, 2003]{Cohen2003}
Cohen, A. (2003).
\newblock {\em Numerical analysis of wavelet methods}, volume~32 of {\em
  Studies in Mathematics and its Applications}.
\newblock North-Holland Publishing Co., Amsterdam.

\bibitem[Courty et~al., 2017]{courty2017}
Courty, N., Flamary, R., Habrard, A., and Rakotomamonjy, A. (2017).
\newblock Joint distribution optimal transportation for domain adaptation.
\newblock In {\em Proceedings of the 31st International Conference on Neural
  Information Processing Systems}, NIPS'17, page 3733–3742. Curran Associates
  Inc.

\bibitem[Courty et~al., 2016]{courty2016optimal}
Courty, N., Flamary, R., Tuia, D., and Rakotomamonjy, A. (2016).
\newblock Optimal transport for domain adaptation.
\newblock {\em IEEE transactions on pattern analysis and machine intelligence},
  39(9):1853--1865.

\bibitem[Cuturi, 2013]{cuturi2013sinkhorn}
Cuturi, M. (2013).
\newblock Sinkhorn distances: Lightspeed computation of optimal transport.
\newblock {\em Advances in neural information processing systems},
  26:2292--2300.

\bibitem[Cuturi and Doucet, 2014]{cuturi2014fast}
Cuturi, M. and Doucet, A. (2014).
\newblock Fast computation of wasserstein barycenters.
\newblock In {\em International conference on machine learning}, pages
  685--693. PMLR.

\bibitem[Damodaran et~al., 2018]{damodaran2018deepjdot}
Damodaran, B.~B., Kellenberger, B., Flamary, R., Tuia, D., and Courty, N.
  (2018).
\newblock Deepjdot: Deep joint distribution optimal transport for unsupervised
  domain adaptation.
\newblock In {\em Proceedings of the European Conference on Computer Vision
  (ECCV)}, pages 447--463.

\bibitem[Deb et~al., 2021]{deb2021efficiency}
Deb, N., Bhattacharya, B.~B., and Sen, B. (2021).
\newblock Efficiency lower bounds for distribution-free hotelling-type
  two-sample tests based on optimal transport.
\newblock {\em arXiv preprint arXiv:2104.01986}.

\bibitem[Deb et~al., 2020]{deb2020measuring}
Deb, N., Ghosal, P., and Sen, B. (2020).
\newblock Measuring association on topological spaces using kernels and
  geometric graphs.
\newblock {\em arXiv preprint arXiv:2010.01768}.

\bibitem[Deb and Sen, 2021]{deb2021multivariate}
Deb, N. and Sen, B. (2021).
\newblock Multivariate rank-based distribution-free nonparametric testing using
  measure transportation.
\newblock {\em Journal of the American Statistical Association},
  (just-accepted):1--45.

\bibitem[Del~Barrio and Loubes, 2019]{del2019}
Del~Barrio, E. and Loubes, J.-M. (2019).
\newblock Central limit theorems for empirical transportation cost in general
  dimension.
\newblock {\em The Annals of Probability}, 47(2):926--951.

\bibitem[Dereich et~al., 2013]{Dereich2013}
Dereich, S., Scheutzow, M., and Schottstedt, R. (2013).
\newblock Constructive quantization: approximation by empirical measures.
\newblock {\em Ann. Inst. Henri Poincar\'{e} Probab. Stat.}, 49(4):1183--1203.

\bibitem[Donoho et~al., 1996]{Donoho1996}
Donoho, D.~L., Johnstone, I.~M., Kerkyacharian, G., and Picard, D. (1996).
\newblock Density estimation by wavelet thresholding.
\newblock {\em Ann. Statist.}, 24(2):508--539.

\bibitem[Drouet~Mari and Kotz, 2001]{Drouet2001}
Drouet~Mari, D. and Kotz, S. (2001).
\newblock {\em Correlation and dependence}.
\newblock Imperial College Press, London; distributed by World Scientific
  Publishing Co., Inc., River Edge, NJ.

\bibitem[Einmahl and Mason, 2000]{Einmahl2000}
Einmahl, U. and Mason, D.~M. (2000).
\newblock An empirical process approach to the uniform consistency of
  kernel-type function estimators.
\newblock {\em J. Theoret. Probab.}, 13(1):1--37.

\bibitem[Ekeland et~al., 2010]{ekeland2010optimal}
Ekeland, I., Galichon, A., and Henry, M. (2010).
\newblock Optimal transportation and the falsifiability of incompletely
  specified economic models.
\newblock {\em Economic Theory}, 42(2):355--374.

\bibitem[El~Moselhy and Marzouk, 2012]{el2012bayesian}
El~Moselhy, T.~A. and Marzouk, Y.~M. (2012).
\newblock Bayesian inference with optimal maps.
\newblock {\em Journal of Computational Physics}, 231(23):7815--7850.

\bibitem[Ferradans et~al., 2014]{ferradans2014regularized}
Ferradans, S., Papadakis, N., Peyr{\'e}, G., and Aujol, J.-F. (2014).
\newblock Regularized discrete optimal transport.
\newblock {\em SIAM Journal on Imaging Sciences}, 7(3):1853--1882.

\bibitem[Folland, 1999]{folland1999real}
Folland, G.~B. (1999).
\newblock {\em Real analysis: modern techniques and their applications},
  volume~40.
\newblock John Wiley \& Sons.

\bibitem[Forrow et~al., 2019]{forrow2019statistical}
Forrow, A., H{\"u}tter, J.-C., Nitzan, M., Rigollet, P., Schiebinger, G., and
  Weed, J. (2019).
\newblock Statistical optimal transport via factored couplings.
\newblock In {\em The 22nd International Conference on Artificial Intelligence
  and Statistics}, pages 2454--2465. PMLR.

\bibitem[Fournier and Guillin, 2015]{Fournier2015}
Fournier, N. and Guillin, A. (2015).
\newblock On the rate of convergence in {W}asserstein distance of the empirical
  measure.
\newblock {\em Probab. Theory Related Fields}, 162(3-4):707--738.

\bibitem[Friesz and Fernandez, 1979]{friesz1979model}
Friesz, T.~L. and Fernandez, J.~E. (1979).
\newblock A model of optimal transport maintenance with demand responsiveness.
\newblock {\em Transportation Research Part B: Methodological}, 13(4):317--339.

\bibitem[Frogner et~al., 2015]{Frogner2015}
Frogner, C., Zhang, C., Mobahi, H., Araya-Polo, M., and Poggio, T. (2015).
\newblock Learning with a {W}asserstein loss.
\newblock In {\em Proceedings of the 28th International Conference on Neural
  Information Processing Systems - Volume 2}, NIPS'15, page 2053–2061,
  Cambridge, MA, USA. MIT Press.

\bibitem[Fukumizu et~al., 2008]{fukumizu2008characteristic}
Fukumizu, K., Sriperumbudur, B.~K., Gretton, A., and Sch{\"o}lkopf, B. (2008).
\newblock Characteristic kernels on groups and semigroups.
\newblock In {\em NIPS}, pages 473--480.

\bibitem[Gabow and Tarjan, 1989]{gabow1989faster}
Gabow, H.~N. and Tarjan, R.~E. (1989).
\newblock Faster scaling algorithms for network problems.
\newblock {\em SIAM Journal on Computing}, 18(5):1013--1036.

\bibitem[Galichon, 2016]{galichon2016optimal}
Galichon, A. (2016).
\newblock {\em Optimal transport methods in economics}.
\newblock Princeton University Press.

\bibitem[Ghosal and Sen, 2019]{ghosal2019multivariate}
Ghosal, P. and Sen, B. (2019).
\newblock Multivariate ranks and quantiles using optimal transportation and
  applications to goodness-of-fit testing.
\newblock {\em arXiv preprint arXiv:1905.05340}.

\bibitem[Gigli, 2011]{Gigli2011}
Gigli, N. (2011).
\newblock On {H}\"{o}lder continuity-in-time of the optimal transport map
  towards measures along a curve.
\newblock {\em Proc. Edinb. Math. Soc. (2)}, 54(2):401--409.

\bibitem[Gin\'{e} and Nickl, 2008]{Gine2008}
Gin\'{e}, E. and Nickl, R. (2008).
\newblock Uniform central limit theorems for kernel density estimators.
\newblock {\em Probab. Theory Related Fields}, 141(3-4):333--387.

\bibitem[Gin\'{e} and Nickl, 2016]{Gine2016}
Gin\'{e}, E. and Nickl, R. (2016).
\newblock {\em Mathematical foundations of infinite-dimensional statistical
  models}.
\newblock Cambridge Series in Statistical and Probabilistic Mathematics, [40].
  Cambridge University Press, New York.

\bibitem[Glaunes et~al., 2004]{glaunes2004diffeomorphic}
Glaunes, J., Trouv{\'e}, A., and Younes, L. (2004).
\newblock Diffeomorphic matching of distributions: A new approach for
  unlabelled point-sets and sub-manifolds matching.
\newblock In {\em Proceedings of the 2004 IEEE Computer Society Conference on
  Computer Vision and Pattern Recognition, 2004. CVPR 2004.}, volume~2, pages
  II--II. IEEE.

\bibitem[Goodfellow et~al., 2014]{goodfellow2014generative}
Goodfellow, I.~J., Pouget-Abadie, J., Mirza, M., Xu, B., Warde-Farley, D.,
  Ozair, S., Courville, A., and Bengio, Y. (2014).
\newblock Generative adversarial networks.
\newblock {\em Advances in Neural Information Processing Systems 27}, pages
  2672--2680.

\bibitem[Gopalan et~al., 2011]{gopalan2011domain}
Gopalan, R., Li, R., and Chellappa, R. (2011).
\newblock Domain adaptation for object recognition: An unsupervised approach.
\newblock In {\em 2011 International Conference on Computer Vision}, pages
  999--1006. IEEE.

\bibitem[Gordaliza et~al., 2019]{gordaliza2019obtaining}
Gordaliza, P., Del~Barrio, E., Fabrice, G., and Loubes, J.-M. (2019).
\newblock Obtaining fairness using optimal transport theory.
\newblock In {\em International Conference on Machine Learning}, pages
  2357--2365. PMLR.

\bibitem[Gretton et~al., 2005a]{gretton2005measuring}
Gretton, A., Bousquet, O., Smola, A., and Sch{\"o}lkopf, B. (2005a).
\newblock Measuring statistical dependence with hilbert-schmidt norms.
\newblock In {\em International conference on algorithmic learning theory},
  pages 63--77. Springer.

\bibitem[Gretton et~al., 2007]{gretton2007kernel}
Gretton, A., Fukumizu, K., Teo, C.~H., Song, L., Sch{\"o}lkopf, B., Smola,
  A.~J., et~al. (2007).
\newblock A kernel statistical test of independence.
\newblock In {\em Nips}, volume~20, pages 585--592. Citeseer.

\bibitem[Gretton et~al., 2005b]{GrettonKernelMeasInd05}
Gretton, A., Herbrich, R., Smola, A., Bousquet, O., and Sch{\"o}lkopf, B.
  (2005b).
\newblock Kernel methods for measuring independence.
\newblock {\em J. Mach. Learn. Res.}, 6:2075--2129 (electronic).

\bibitem[Gunsilius and Schennach, 2019]{gunsilius2019independent}
Gunsilius, F. and Schennach, S.~M. (2019).
\newblock Independent nonlinear component analysis.
\newblock Technical report, cemmap working paper.

\bibitem[Gunsilius, 2021]{gunsilius2018convergence}
Gunsilius, F.~F. (2021).
\newblock On the convergence rate of potentials of brenier maps.
\newblock {\em Econometric Theory}, pages 1--37.

\bibitem[Guntuboyina and Sen, 2012]{guntuboyina2012}
Guntuboyina, A. and Sen, B. (2012).
\newblock L1 covering numbers for uniformly bounded convex functions.
\newblock In {\em Conference on Learning Theory}, pages 12--1. JMLR Workshop
  and Conference Proceedings.

\bibitem[Hallin et~al., 2021]{hallin2021distribution}
Hallin, M., Del~Barrio, E., Cuesta-Albertos, J., and Matr{\'a}n, C. (2021).
\newblock Distribution and quantile functions, ranks and signs in dimension d:
  A measure transportation approach.
\newblock {\em The Annals of Statistics}, 49(2):1139--1165.

\bibitem[Hansen, 2008]{Hansen08}
Hansen, B.~E. (2008).
\newblock Uniform convergence rates for kernel estimation with dependent data.
\newblock {\em Econometric Theory}, 24(3):726--748.

\bibitem[H\"{a}rdle et~al., 1998]{Hardle1998}
H\"{a}rdle, W., Kerkyacharian, G., Picard, D., and Tsybakov, A. (1998).
\newblock {\em Wavelets, approximation, and statistical applications}, volume
  129 of {\em Lecture Notes in Statistics}.
\newblock Springer-Verlag, New York.

\bibitem[Heller et~al., 2013]{Heller2013}
Heller, R., Heller, Y., and Gorfine, M. (2013).
\newblock A consistent multivariate test of association based on ranks of
  distances.
\newblock {\em Biometrika}, 100(2):503--510.

\bibitem[Hiriart-Urruty and Lemar\'{e}chal, 1993]{hiriart1993}
Hiriart-Urruty, J.-B. and Lemar\'{e}chal, C. (1993).
\newblock {\em Convex analysis and minimization algorithms. {II}}, volume 306
  of {\em Grundlehren der Mathematischen Wissenschaften [Fundamental Principles
  of Mathematical Sciences]}.
\newblock Springer-Verlag, Berlin.
\newblock Advanced theory and bundle methods.

\bibitem[H{\"u}tter and Rigollet, 2021]{hutter2021minimax}
H{\"u}tter, J.-C. and Rigollet, P. (2021).
\newblock Minimax estimation of smooth optimal transport maps.
\newblock {\em The Annals of Statistics}, 49(2):1166--1194.

\bibitem[Jonker and Volgenant, 1987]{jonker1987shortest}
Jonker, R. and Volgenant, A. (1987).
\newblock A shortest augmenting path algorithm for dense and sparse linear
  assignment problems.
\newblock {\em Computing}, 38(4):325--340.

\bibitem[Josse and Holmes, 2016]{Josse2016}
Josse, J. and Holmes, S. (2016).
\newblock Measuring multivariate association and beyond.
\newblock {\em Stat. Surv.}, 10:132--167.

\bibitem[Kandasamy et~al., 2018]{Kandaswamy2018}
Kandasamy, K., Neiswanger, W., Schneider, J., P\'{o}czos, B., and Xing, E.~P.
  (2018).
\newblock Neural architecture search with bayesian optimisation and optimal
  transport.
\newblock In {\em Proceedings of the 32nd International Conference on Neural
  Information Processing Systems}, NIPS'18, page 2020–2029, Red Hook, NY,
  USA. Curran Associates Inc.

\bibitem[Kantorovich, 2004]{Kantorovic1948}
Kantorovich, L.~V. (2004).
\newblock On a problem of {M}onge.
\newblock {\em Zap. Nauchn. Sem. S.-Peterburg. Otdel. Mat. Inst. Steklov.
  (POMI)}, 312(Teor. Predst. Din. Sist. Komb. i Algoritm. Metody. 11):15--16.

\bibitem[Kantorovitch, 1942]{Kantorovic1942}
Kantorovitch, L. (1942).
\newblock On the translocation of masses.
\newblock {\em C. R. (Doklady) Acad. Sci. URSS (N.S.)}, 37:199--201.

\bibitem[Kerkyacharian and Picard, 1992]{Kerk1992}
Kerkyacharian, G. and Picard, D. (1992).
\newblock Density estimation in {B}esov spaces.
\newblock {\em Statist. Probab. Lett.}, 13(1):15--24.

\bibitem[Kim et~al., 2020]{kim2020minimax}
Kim, I., Balakrishnan, S., and Wasserman, L. (2020).
\newblock Minimax optimality of permutation tests.
\newblock {\em arXiv preprint arXiv:2003.13208}.

\bibitem[Kim et~al., 2013]{kim2013efficient}
Kim, S., Ma, R., Mesa, D., and Coleman, T.~P. (2013).
\newblock Efficient bayesian inference methods via convex optimization and
  optimal transport.
\newblock In {\em 2013 IEEE International Symposium on Information Theory},
  pages 2259--2263. IEEE.

\bibitem[Kingma and Welling, 2014]{kingma2013auto}
Kingma, D.~P. and Welling, M. (2014).
\newblock {Auto-Encoding Variational Bayes}.
\newblock In {\em 2nd International Conference on Learning Representations,
  {ICLR} 2014, Banff, AB, Canada, April 14-16, 2014, Conference Track
  Proceedings}.

\bibitem[Koml\'{o}s et~al., 1976]{Komlos1976}
Koml\'{o}s, J., Major, P., and Tusn\'{a}dy, G. (1976).
\newblock An approximation of partial sums of independent {RV}'s, and the
  sample {DF}. {II}.
\newblock {\em Z. Wahrscheinlichkeitstheorie und Verw. Gebiete}, 34(1):33--58.

\bibitem[Kur et~al., 2020]{kur2020convex}
Kur, G., Gao, F., Guntuboyina, A., and Sen, B. (2020).
\newblock Convex regression in multidimensions: Suboptimality of least squares
  estimators.
\newblock {\em arXiv preprint arXiv:2006.02044}.

\bibitem[Laplace, 1820]{laplace1820theorie}
Laplace, P.~S. (1820).
\newblock {\em Th{\'e}orie analytique des probabilit{\'e}s}.
\newblock Courcier.

\bibitem[Le~Gouic and Loubes, 2017]{le2017existence}
Le~Gouic, T. and Loubes, J.-M. (2017).
\newblock Existence and consistency of wasserstein barycenters.
\newblock {\em Probability Theory and Related Fields}, 168(3):901--917.

\bibitem[Li and Yuan, 2019]{li2019optimality}
Li, T. and Yuan, M. (2019).
\newblock On the optimality of gaussian kernel based nonparametric tests
  against smooth alternatives.
\newblock {\em arXiv preprint arXiv:1909.03302}.

\bibitem[Li and Nochetto, 2020]{li2020quantitative}
Li, W. and Nochetto, R.~H. (2020).
\newblock Quantitative stability and error estimates for optimal transport
  plans.
\newblock {\em IMA Journal of Numerical Analysis}.

\bibitem[Li et~al., 2015]{li2015generative}
Li, Y., Swersky, K., and Zemel, R. (2015).
\newblock Generative moment matching networks.
\newblock In {\em International Conference on Machine Learning}, pages
  1718--1727. PMLR.

\bibitem[Luise et~al., 2018]{Luise2018}
Luise, G., Rudi, A., Pontil, M., and Ciliberto, C. (2018).
\newblock Differential properties of sinkhorn approximation for learning with
  {W}asserstein distance.
\newblock In {\em Proceedings of the 32nd International Conference on Neural
  Information Processing Systems}, NIPS'18, page 5864–5874, Red Hook, NY,
  USA. Curran Associates Inc.

\bibitem[Masarotto et~al., 2019]{masarotto2019procrustes}
Masarotto, V., Panaretos, V.~M., and Zemel, Y. (2019).
\newblock Procrustes metrics on covariance operators and optimal transportation
  of gaussian processes.
\newblock {\em Sankhy\={a} A}, 81(1):172--213.

\bibitem[Mason, 2012]{Mason2012}
Mason, D.~M. (2012).
\newblock Proving consistency of non-standard kernel estimators.
\newblock {\em Stat. Inference Stoch. Process.}, 15(2):151--176.

\bibitem[McCann, 1995]{McCann95}
McCann, R.~J. (1995).
\newblock Existence and uniqueness of monotone measure-preserving maps.
\newblock {\em Duke Math. J.}, 80(2):309--323.

\bibitem[Mena and Niles-Weed, 2019]{Gonzalo2019}
Mena, G. and Niles-Weed, J. (2019).
\newblock Statistical bounds for entropic optimal transport: Sample complexity
  and the central limit theorem.
\newblock {\em Advances in Neural Information Processing Systems}, 32.

\bibitem[M{\'e}rigot, 2011]{Merigot11}
M{\'e}rigot, Q. (2011).
\newblock A multiscale approach to optimal transport.
\newblock In {\em Computer Graphics Forum}, volume~30, pages 1583--1592. Wiley
  Online Library.

\bibitem[Merigot and Thibert, 2020]{Merigot2020}
Merigot, Q. and Thibert, B. (2020).
\newblock Optimal transport: discretization and algorithms.
\newblock {\em Handbook of Numerical Analysis 22 -- Geometric PDES (arXiv
  preprint arXiv:2003.00855)}.

\bibitem[Meyer, 1990]{Meyer1990}
Meyer, Y. (1990).
\newblock {\em Ondelettes et op\'{e}rateurs. {I}}.
\newblock Actualit\'{e}s Math\'{e}matiques. [Current Mathematical Topics].
  Hermann, Paris.
\newblock Ondelettes. [Wavelets].

\bibitem[Mohamed and Lakshminarayanan, 2017]{shakir2017}
Mohamed, S. and Lakshminarayanan, B. (2017).
\newblock Learning in implicit generative models.
\newblock In {\em Proceedings of the International Conference in Learning
  Representations}.

\bibitem[Monge, 1781]{Monge1781}
Monge, G. (1781).
\newblock M\'{e}moire sur la th\'{e}orie des d\'{e}blais et des remblais.
\newblock {\em Histoire de l'Acad\'{e}mie Royale des Sciences de Paris}, pages
  666--704.

\bibitem[M\"{u}ller, 1997]{Alfred1997}
M\"{u}ller, A. (1997).
\newblock Integral probability metrics and their generating classes of
  functions.
\newblock {\em Adv. in Appl. Probab.}, 29(2):429--443.

\bibitem[Nadaraja, 1965]{Nadaraja1965}
Nadaraja, E.~A. (1965).
\newblock On non-parametric estimates of density functions and regression.
\newblock {\em Teor. Verojatnost. i Primenen.}, 10:199--203.

\bibitem[Niles-Weed and Rigollet, 2019]{niles2019estimation}
Niles-Weed, J. and Rigollet, P. (2019).
\newblock Estimation of {W}asserstein distances in the spiked transport model.
\newblock {\em arXiv preprint arXiv:1909.07513}.

\bibitem[Oudjane and Musso, 2005]{oudjane2005sup}
Oudjane, N. and Musso, C. (2005).
\newblock L2-density estimation with negative kernels.
\newblock In {\em ISPA 2005. Proceedings of the 4th International Symposium on
  Image and Signal Processing and Analysis, 2005.}, pages 34--39. IEEE.

\bibitem[Papadakis et~al., 2014]{papadakis2014optimal}
Papadakis, N., Peyr{\'e}, G., and Oudet, E. (2014).
\newblock Optimal transport with proximal splitting.
\newblock {\em SIAM Journal on Imaging Sciences}, 7(1):212--238.

\bibitem[Parzen, 1962]{Parzen1962}
Parzen, E. (1962).
\newblock On estimation of a probability density function and mode.
\newblock {\em Ann. Math. Statist.}, 33:1065--1076.

\bibitem[Peyr{\'e} et~al., 2019]{peyre2019computational}
Peyr{\'e}, G., Cuturi, M., et~al. (2019).
\newblock Computational optimal transport: With applications to data science.
\newblock {\em Foundations and Trends{\textregistered} in Machine Learning},
  11(5-6):355--607.

\bibitem[Rachev, 1985]{rachev1985monge}
Rachev, S.~T. (1985).
\newblock The monge--kantorovich mass transference problem and its stochastic
  applications.
\newblock {\em Theory of Probability \& Its Applications}, 29(4):647--676.

\bibitem[Radulovi\'{c} and Wegkamp, 2000]{Radulovic2000}
Radulovi\'{c}, D. and Wegkamp, M. (2000).
\newblock Weak convergence of smoothed empirical processes: beyond {D}onsker
  classes.
\newblock In {\em High dimensional probability, {II} ({S}eattle, {WA}, 1999)},
  volume~47 of {\em Progr. Probab.}, pages 89--105. Birkh\"{a}user Boston,
  Boston, MA.

\bibitem[Ramdas et~al., 2017]{Ramdas2017wasserstein}
Ramdas, A., Trillos, N.~G., and Cuturi, M. (2017).
\newblock On {W}asserstein two-sample testing and related families of
  nonparametric tests.
\newblock {\em Entropy}, 19(2):47.

\bibitem[Reich, 2011]{reich2011dynamical}
Reich, S. (2011).
\newblock A dynamical systems framework for intermittent data assimilation.
\newblock {\em BIT Numerical Mathematics}, 51(1):235--249.

\bibitem[Reich, 2013]{reich2013nonparametric}
Reich, S. (2013).
\newblock A nonparametric ensemble transform method for bayesian inference.
\newblock {\em SIAM Journal on Scientific Computing}, 35(4):A2013--A2024.

\bibitem[Rippl et~al., 2016]{Rippl2016}
Rippl, T., Munk, A., and Sturm, A. (2016).
\newblock Limit laws of the empirical {W}asserstein distance: {G}aussian
  distributions.
\newblock {\em J. Multivariate Anal.}, 151:90--109.

\bibitem[Salimans et~al., 2018]{salimans2018improving}
Salimans, T., Zhang, H., Radford, A., and Metaxas, D. (2018).
\newblock Improving {GAN}s using optimal transport.
\newblock {\em arXiv preprint arXiv:1803.05573}.

\bibitem[Santambrogio, 2015]{santambrogio2015optimal}
Santambrogio, F. (2015).
\newblock Optimal transport for applied mathematicians.
\newblock {\em Birk{\"a}user, NY}, 55(58-63):94.

\bibitem[Seguy et~al., 2018]{seguy2018large}
Seguy, V., Damodaran, B.~B., Flamary, R., Courty, N., Rolet, A., and Blondel,
  M. (2018).
\newblock Large-scale optimal transport and mapping estimation.
\newblock In {\em Proceedings of the International Conference in Learning
  Representations}.

\bibitem[Shi et~al., 2020a]{shi2020distribution}
Shi, H., Drton, M., and Han, F. (2020a).
\newblock Distribution-free consistent independence tests via center-outward
  ranks and signs.
\newblock {\em Journal of the American Statistical Association}, pages 1--16.

\bibitem[Shi et~al., 2020b]{Shi2020}
Shi, H., Hallin, M., Drton, M., and Han, F. (2020b).
\newblock Rate-optimality of consistent distribution-free tests of independence
  based on center-outward ranks and signs.
\newblock {\em arXiv preprint arXiv:2007.02186}.

\bibitem[Silverman, 1978]{Silverman1978}
Silverman, B.~W. (1978).
\newblock Weak and strong uniform consistency of the kernel estimate of a
  density and its derivatives.
\newblock {\em Ann. Statist.}, 6(1):177--184.

\bibitem[Sommerfeld and Munk, 2018]{Sommerfeld2018}
Sommerfeld, M. and Munk, A. (2018).
\newblock Inference for empirical {W}asserstein distances on finite spaces.
\newblock {\em J. R. Stat. Soc. Ser. B. Stat. Methodol.}, 80(1):219--238.

\bibitem[Sriperumbudur et~al., 2012]{Sriperumbudur2012}
Sriperumbudur, B.~K., Fukumizu, K., Gretton, A., Sch\"{o}lkopf, B., and
  Lanckriet, G. R.~G. (2012).
\newblock On the empirical estimation of integral probability metrics.
\newblock {\em Electron. J. Stat.}, 6:1550--1599.

\bibitem[Sriperumbudur et~al., 2010]{Sriperumbudur2010}
Sriperumbudur, B.~K., Gretton, A., Fukumizu, K., Sch\"{o}lkopf, B., and
  Lanckriet, G. R.~G. (2010).
\newblock Hilbert space embeddings and metrics on probability measures.
\newblock {\em J. Mach. Learn. Res.}, 11:1517--1561.

\bibitem[Srivastava et~al., 2015]{srivastava2015wasp}
Srivastava, S., Cevher, V., Dinh, Q., and Dunson, D. (2015).
\newblock Wasp: Scalable bayes via barycenters of subset posteriors.
\newblock In {\em Artificial Intelligence and Statistics}, pages 912--920.
  PMLR.

\bibitem[Srivastava et~al., 2018]{srivastava2018scalable}
Srivastava, S., Li, C., and Dunson, D.~B. (2018).
\newblock Scalable bayes via barycenter in wasserstein space.
\newblock {\em The Journal of Machine Learning Research}, 19(1):312--346.

\bibitem[Su et~al., 2015]{su2015optimal}
Su, Z., Wang, Y., Shi, R., Zeng, W., Sun, J., Luo, F., and Gu, X. (2015).
\newblock Optimal mass transport for shape matching and comparison.
\newblock {\em IEEE transactions on pattern analysis and machine intelligence},
  37(11):2246--2259.

\bibitem[Sz{\'e}kely et~al., 2007]{SzekelyCorrDist07}
Sz{\'e}kely, G.~J., Rizzo, M.~L., and Bakirov, N.~K. (2007).
\newblock Measuring and testing dependence by correlation of distances.
\newblock {\em Ann. Statist.}, 35(6):2769--2794.

\bibitem[Talagrand, 1994]{Talagrand1994}
Talagrand, M. (1994).
\newblock The transportation cost from the uniform measure to the empirical
  measure in dimension {$\ge 3$}.
\newblock {\em Ann. Probab.}, 22(2):919--959.

\bibitem[Triebel, 2006]{Triebel2006}
Triebel, H. (2006).
\newblock {\em Theory of function spaces. {III}}, volume 100 of {\em Monographs
  in Mathematics}.
\newblock Birkh\"{a}user Verlag, Basel.

\bibitem[Tsybakov, 2009]{Tsybakov2009}
Tsybakov, A.~B. (2009).
\newblock {\em Introduction to nonparametric estimation}.
\newblock Springer Series in Statistics. Springer, New York.
\newblock Revised and extended from the 2004 French original, Translated by
  Vladimir Zaiats.

\bibitem[van~de Geer, 2000]{geer2000}
van~de Geer, S.~A. (2000).
\newblock {\em Applications of empirical process theory}, volume~6 of {\em
  Cambridge Series in Statistical and Probabilistic Mathematics}.
\newblock Cambridge University Press, Cambridge.

\bibitem[van~der Vaart, 1998]{vdV98}
van~der Vaart, A.~W. (1998).
\newblock {\em Asymptotic statistics}, volume~3 of {\em Cambridge Series in
  Statistical and Probabilistic Mathematics}.
\newblock Cambridge University Press, Cambridge.

\bibitem[van~der Vaart and Wellner, 1996]{vaart2013}
van~der Vaart, A.~W. and Wellner, J.~A. (1996).
\newblock {\em Weak convergence and empirical processes}.
\newblock Springer Series in Statistics. Springer-Verlag, New York.
\newblock With applications to statistics.

\bibitem[Villani, 2003]{V03}
Villani, C. (2003).
\newblock {\em Topics in optimal transportation}, volume~58 of {\em Graduate
  Studies in Mathematics}.
\newblock American Mathematical Society, Providence, RI.

\bibitem[Villani, 2009]{V09}
Villani, C. (2009).
\newblock {\em Optimal transport}, volume 338 of {\em Grundlehren der
  Mathematischen Wissenschaften [Fundamental Principles of Mathematical
  Sciences]}.
\newblock Springer-Verlag, Berlin.
\newblock Old and new.

\bibitem[Wainwright, 2019]{Wainwright2019}
Wainwright, M.~J. (2019).
\newblock {\em High-dimensional statistics}, volume~48 of {\em Cambridge Series
  in Statistical and Probabilistic Mathematics}.
\newblock Cambridge University Press, Cambridge.
\newblock A non-asymptotic viewpoint.

\bibitem[Walter, 1992]{Walter1992}
Walter, G.~G. (1992).
\newblock Approximation of the delta function by wavelets.
\newblock {\em J. Approx. Theory}, 71(3):329--343.

\bibitem[Weed and Berthet, 2019]{weed2019estimation}
Weed, J. and Berthet, Q. (2019).
\newblock Estimation of smooth densities in {W}asserstein distance.
\newblock In {\em Conference on Learning Theory}, pages 3118--3119. PMLR.

\bibitem[Zemel and Panaretos, 2019]{zemel2019frechet}
Zemel, Y. and Panaretos, V.~M. (2019).
\newblock Fr\'{e}chet means and {P}rocrustes analysis in {W}asserstein space.
\newblock {\em Bernoulli}, 25(2):932--976.

\end{thebibliography}

	\appendix
	
	\section{Appendix}
	This section is devoted to proving our main results and is organized as follows: In~\cref{sec:mainrespf}, we present the proofs of results from~\cref{sec:mainres} and in~\cref{sec:Appf}, we present the proofs from~\cref{sec:App}. Throughout this section, we will use the $\lesssim$ sign to hide constants that are free of $m,n$.
	\subsection{Proofs from~\cref{sec:mainres}}\label{sec:mainrespf}
	\begin{proof}[Proof of~\cref{thm:newubd}]
		We begin the proof by observing that $\phmn^*(\cdot)$ is convex and finite on $\mbox{supp}(\nu)$, and hence differentiable $\nu$ almost everywhere (a.e.). Further by~\cref{lem:gradual}, we also have:
		\begin{equation}\label{eq:newubd1}
		\nabla\phmn^*(\Tmn(\mx))=\mx \qquad \mbox{$\mu$-a.e.~$\mx$.}
		\end{equation} 
		Fix any arbitrary $\gamma\in\tgm$ and suppose that $\gamma(\my|\mx)$ denotes the conditional distribution of $\my$ given $\mx$ under $\gamma$. Define,
		$$D_1:= \int  \phmn^*(\my)\,d\tnu(\my) - \int  \phmn^*(\my) d \bar{\nu}_m(\my).$$
		As $\gamma$ has marginals $\tmu$ and $\tnu$, we have:
		\begin{align}\label{eq:newubd2}
		D_1=\int_{\mx,\my} \phmn^*(\my)\,d\gamma(\my|\mx)\,d\tmu(\mx)-\int_{\mx} \phmn^*(\Tmn(\mx))\,d\tmu(\mx).
		\end{align}
		Next, by applying the conditional version of Jensen's inequality,
		\begin{align}\label{eq:newubd3}
		\int_{\mx} \left(\int_{\my} \phmn^*(\my)\,d\gamma(\my|\mx)\right)\,d\tmu(\mx)&\geq \int_{\mx}\phmn^*\left(\int_{\my} \my\,d\gamma(\my|\mx)\right)\,d\tmu(\mx)\nonumber \\&=\int_{\mx}\phmn^*(\ttg(\mx))\,d\tmu(\mx).
		\end{align}
		Using~\eqref{eq:newubd3} with~\eqref{eq:newubd2} yields,
		\begin{align}\label{eq:newubd4}
		D_1
		&\geq \int [\phmn^*(\ttg(\mx)) - \phmn^*(\Tmn(\mx))] d \tmu(\mx)  \nonumber \\ 
		& \overset{(a)}{\ge} \int \left\{\nabla \phmn^*(\Tmn(\mx))^\top (\ttg(\mx) - \Tmn(\mx)) + \frac{1}{2L} \|\ttg(\mx)-\Tmn(\mx)\|^2 \right\} \,d {\tmu}(\mx)  \nonumber \\
		& \overset{(b)}{=} \underbrace{\int \mx^\top (\ttg(\mx) - \Tmn(\mx)) \,d {\tmu}(\mx)}_{D_2} + \frac{1}{2L} \int  \|\ttg(\mx)-\Tmn(\mx)\|^2 \, d {\tmu}(\mx).
		\end{align}
		Here (a) follows from the strong convexity of $\phmn^*(\cdot)$ with parameter $(1/L)$ (see~\cref{lem:assm}) and (b) follows from~\eqref{eq:newubd1}.
		
		Next, we will simplify the term $D_2$. Towards this direction, observe that for every $\gamma\in\tgm$,
		\begin{align}\label{eq:newubd5}
		W_2^2(\tmu,\tnu)&=\int \lVert \mx-\my\rVert^2\,d\gamma(\mx,\my)\nonumber \\ &=\int \lVert \mx\rVert^2\,d\tmu(\mx)+\int \lVert \my\rVert^2\,d\tnu(\my)-2\int_{\mx} \left(\mx^{\top}\int_{\my} \my\,d\gamma(\my|\mx)\right)\,d\tmu(\mx)\nonumber \\&=\int \lVert \mx\rVert^2\,d\tmu(\mx)+\int \lVert \my\rVert^2\,d\tnu(\my)-2\int_{\mx} \mx^{\top}\ttg(\mx)\,d\tmu(\mx).
		\end{align}
		Also, as $\Tmn$ is the gradient of a convex function, it is also an OT map from $\tmu$ to $\bar{\nu}_m$ (see~\cite[Section 1.2]{Ambrosio2013}), we have:
		\begin{align}\label{eq:newubd6}
		W_2^2(\tmu,\bar{\nu}_m)&=\int \lVert \mx-\Tmn(\mx)\rVert^2\,d\tmu(\mx)\nonumber \\ &=\int \lVert \mx\rVert^2\,d\tmu(\mx)+\int \lVert \my\rVert^2\,d\bar{\nu}_m(\my)-2\int_{\mx} \mx^{\top}\Tmn(\mx)\,d\tmu(\mx).
		\end{align}
		Now~\eqref{eq:newubd5}~and~\eqref{eq:newubd6} imply
		\begin{align}\label{eq:newubd7}
		D_2=\frac{1}{2}\big(W_2^2(\tmu,\bar{\nu}_m)-W_2^2(\tmu,\tnu)\big)+\frac{1}{2}\int \lVert \my\rVert^2\,d(\tnu-\bar{\nu}_m)(\my).
		\end{align}
		Finally by combining~\eqref{eq:newubd7}~and~\eqref{eq:newubd4}, we get:
		\begin{align}\label{eq:newubd8}&\;\;\;\;\frac{1}{2L}\int \lVert \ttg(\mx)-\Tmn(\mx)\rVert^2\,d\tmu(\mx)\nonumber \\ &\leq \frac{1}{2}\big(W_2^2(\tmu,\tnu)-W_2^2(\tmu,\bar{\nu}_m)\big)+\int (\phmn^*(\my)-(1/2)\lVert \my\rVert^2)\,d(\tnu-\bar{\nu}_m)(\my).\end{align}
		Now note that the bound on the right hand side of the above display is free of the particular choice of $\gamma\in\tgm$. Therefore, the same bound holds if we take a supremum over $\gamma\in\tgm$ on the left hand side. We will now provide an upper bound for the right hand side of~\eqref{eq:newubd8}. The remainder of the proof proceeds as in the proof of~\cite[Proposition 2]{Gonzalo2019}.
		
		By the dual representation presented in~\eqref{eq:altwass} and~\eqref{eq:dualf}, and the definitions of $\pst(\cdot)$ and $\psw(\cdot)$ in the statement of~\cref{thm:newubd}, we have
		$$\frac{1}{2}W_2^2(\tmu,\tnu)=\frac{1}{2}\int \lVert x\rVert^2\,d\tmu(x)+\frac{1}{2}\int\lVert y\rVert^2\,d\tnu(y)-\cS_{\tmu,\tnu}(\pst),$$
		$$\mbox{and}\qquad \frac{1}{2}W_2^2(\tmu,\bar{\nu}_m)=\frac{1}{2}\int \lVert x\rVert^2\,d\tmu(x)+\frac{1}{2}\int\lVert y\rVert^2\,d\bar{\nu}_m(y)-\cS_{\tmu,\bar{\nu}_m}(\psw).$$
		By subtracting the two equations above, we get:
		\begin{equation}\label{eq:newubd9}
		\frac{1}{2}W_2^2(\tmu,\tnu)-\frac{1}{2}W_2^2(\tmu,\bar{\nu}_m)=\frac{1}{2}\int \lVert y\rVert^2\,d(\tnu-\bar{\nu}_m)-\cS_{\tmu,\tnu}(\pst)+\cS_{\tmu,\bar{\nu}_m}(\psw).
		\end{equation}
		Next, we use~\eqref{eq:dualf} to make the following observations:
		\begin{equation}\label{eq:newubd10}
		\cS_{\tmu,\tnu}(\pst)\leq \cS_{\tmu,\tnu}(\psw),\qquad \cS_{\tmu,\bar{\nu}_m}(\psw)\leq \cS_{\tmu,\bar{\nu}_m}(\pst).
		\end{equation}
		Note that~\eqref{eq:newubd10} immediately yields the following conclusions:
		\begin{equation*}
		\cS_{\tmu,\bar{\nu}_m}(\psw)-\cS_{\tmu,\tnu}(\psw)\leq \cS_{\tmu,\bar{\nu}_m}(\psw)-\cS_{\tmu,\tnu}(\pst),
		\end{equation*}
		and
		\begin{equation*}
		\cS_{\tmu,\bar{\nu}_m}(\psw)-\cS_{\tmu,\tnu}(\pst)\leq \cS_{\tmu,\bar{\nu}_m}(\pst)-\cS_{\tmu,\tnu}(\pst).
		\end{equation*}
		By combining the above two displays, we have:
		\begin{align}\label{eq:newubd11}
		&\;\;\;\;\left|\cS_{\tmu,\bar{\nu}_m}(\psw)-\cS_{\tmu,\tnu}(\pst)\right|\nonumber \\ &\leq \max\left\{\left|\cS_{\tmu,\bar{\nu}_m}(\pst)-\cS_{\tmu,\tnu}(\pst)\right|,\left|\cS_{\tmu,\bar{\nu}_m}(\psw)-\cS_{\tmu,\tnu}(\psw)\right|\right\}.
		\end{align}
		By~\eqref{eq:dualf} and some simple algebra, the following holds:
		$$ \left|\cS_{\tmu,\bar{\nu}_m}(\pst)-\cS_{\tmu,\tnu}(\pst)\right|=\left|\int \pst^*\,d(\bar{\nu}_m-\tnu)\right|.$$
		A similar expression holds for $|\cS_{\tmu,\bar{\nu}_m}(\psw)-\cS_{\tmu,\tnu}(\psw)|$. Using the above observation in~\eqref{eq:newubd11}, we get:
		$$\left|\cS_{\tmu,\bar{\nu}_m}(\psw)-\cS_{\tmu,\tnu}(\pst)\right|\leq \max\left\{\bigg|\int \pst^* \,d(\tnu-\bar{\nu}_m)\bigg|,\bigg|\int \psw^* \,d(\tnu-\bar{\nu}_m)\bigg|\right\}.$$
		Combining the above display with~\eqref{eq:newubd9}, we further have:
		\begin{align}\label{eq:newubd12}
		&\;\;\;\;\left|\frac{1}{2}W_2^2(\tmu,\tnu)-\frac{1}{2}W_2^2(\tmu,\bar{\nu}_m)-\left(\frac{1}{2}\int \lVert y\rVert^2\,d(\tnu-\bar{\nu}_m)\right)\right|\nonumber\\ &\leq \max\left\{\bigg|\int \pst^* \,d(\tnu-\bar{\nu}_m)\bigg|,\bigg|\int \psw^* \,d(\tnu-\bar{\nu}_m)\bigg|\right\}.
		\end{align}
		Combining~\eqref{eq:newubd12}~with~\eqref{eq:newubd8} then completes the proof.
	\end{proof}
	
	\begin{proof}[Proof of~\cref{thm:nsmooth}]				
		First observe that
		$$\limsup_{M\to\infty}\limsup_{m,n\to\infty} \P\left(\Big|\int \phmn^*\,d(\hnu-\bar{\nu}_m)\big|\geq M \left(r_d^{(m,m)}+r_d^{(n,n)}\right)\right)=0$$
		by the weak law of large numbers as $(r_{d}^{(n,n)})^{-1}n^{-1/2}=O(1)$ and $(r_{d}^{(m,m)})^{-1}m^{-1/2}=O(1)$.		
		
		Combining the above observation with~\cref{thm:newubd}, we have:
		\begin{align}\label{eq:nsmooth1}
		&\;\;\;\;\limsup_{M\to\infty}\limsup_{m,n\to\infty}\P\left(\sup\limits_{\gamma\in\tgm} \int \lVert \ttg(\mx)-\Tmn(\mx)\rVert^2 \,d\tmu(\mx)\geq Mr_d^{(m,n)}\right)\nonumber \\ &\leq  \limsup_{M\to\infty}\limsup_{m,n\to\infty}\P\left(\max\left\{\bigg|\int \pst^* \,d(\tnu-\bar{\nu}_m)\bigg|,\bigg|\int \psw^* \,d(\tnu-\bar{\nu}_m)\bigg|\right\}\geq \frac{M}{2} r_d^{(m,n)}\right)\nonumber \\ &\leq \limsup_{M\to\infty}\limsup_{m,n\to\infty}\Bigg[\P\left(\left|\int \pst^* \,d(\tnu-\nu)\right|\geq \frac{M}{2}r_d^{(n,n)}\right)+\P\left(\left|\int \pst^* \,d(\bar{\nu}_m-\nu)\right|\geq \frac{M}{2}r_d^{(m,m)}\right)\nonumber \\&+\P\left(\left|\int \psw^* \,d(\tnu-\nu)\right|\geq \frac{M}{2} r_d^{(n,n)}\right)+\P\left(\left|\int \psw^* \,d(\bar{\nu}_m-\nu)\right|\geq \frac{M}{2} r_d^{(m,m)}\right)\Bigg].
		\end{align}
		In the sequel, we will only discuss how to bound the first term on the right hand side of~\eqref{eq:nsmooth1}. Once that is understood, the other terms can be bounded similarly. Therefore, our focus is on bounding
		\begin{equation}\label{eq:nsmooth2}
		\limsup_{M\to\infty}\limsup_{m,n\to\infty}\P\left(\left|\int \pst^* \,d(\tnu-\nu)\right|\geq \frac{M}{2}r_d^{(n,n)}(\log{(1+\max\{m,n\})})^{t_{d,\alpha}}\right).
		\end{equation}
		For the next part, to simplify notation, let us begin with some notation. Set $\xn:=\mbox{supp}(\nu)$ and $\xmo$ denote the closure of the convex hull of $X_1,\ldots ,X_n$. 
		
		Note that if we replace $\pst(\cdot)$ by $\pst(\cdot)-C$ for some constant $C>0$, then $\pst^*(\cdot)\mapsto\pst^*+C$. However replacing $\pst^*(\cdot)$ by $\pst^*(\cdot)+C$ in~\eqref{eq:nsmooth2} doesn't change its value as $\tnu$ and $\nu$ are both probability measures. Therefore, without loss of generality, we can assume that $\pst(X_1)=0$ for all $m,n$. We will stick to this convention for the rest of the proof. Also note that $\pst(\cdot)$ is only determined at the data points $X_1,\ldots ,X_n$. Without loss of generality, we extend $\pst(\cdot)$ to the whole of $\R^d$ by linear interpolation for any $x\in\xmo$ and setting $\pst(x)=\infty$ for $x\in\xmo^c$.
		
		The proof now proceeds using the following steps:
		
		\textbf{Step I:} There exists a constant $C_1>0$ and  $y_n\in\mbox{supp}(\nu)=\xn$ such that $$|\pst^*(y_n)|\leq \max_{1\leq i\leq m}\lVert X_i\rVert.$$
		\begin{proof}[Proof of step I]
			By Kantorovich duality, there exists $y_n$ such that $$\pst^*(y_n)+\pst(X_1)=\langle X_1,y_n\rangle\quad \implies \quad |\pst^*(y_n)|\leq C_1\lVert X_1\rVert\leq C_1\max_{1\leq i\leq m}\lVert X_i\rVert,$$
			where $C_1:=\sup\{\lVert y\rVert:\ y\in\xn\}$.
		\end{proof}
		
		\textbf{Step II:} There exists a constant $C_2>0$ such that the following holds:
		$$\lVert \pst^*\rVert_{\infty,\xn}\leq C_2\max_{1\leq i\leq n} \lVert X_i\rVert,$$
		where $\lVert\cdot\rVert_{\infty,\xn}$ is the uniform norm on the support of $\nu$.
		\begin{proof}[Proof of step II]
			As $\pst(x)=\infty$ for $x\in\xmo^c$, using~\eqref{eq:lfdual}, we can write $\pst^*(y)=\max_{x\in\xmo} (\langle x,y\rangle-\pst(x))$ for all $y\in\xn$. For any $y_0\in\xn$, let $x_0\in\xmo$ be such that $\pst^*(y_0)=\langle x_0,y_0\rangle-\pst(x_0)$. Then, for any $y\in\X$, we have:
			\begin{align*}
			&\begin{cases} \pst^*(y_0)=\langle x_0,y_0\rangle-\pst(x_0)\\ \pst^*(y)\geq \langle x_0,y\rangle-\pst(x_0)\end{cases}\\ \implies &|\pst^*(y_0)-\pst^*(y)|\leq |\langle x_0,y_0-y\rangle|\leq \left(\max_{1\leq i\leq m}\lVert X_i\rVert\right)\lVert y_0-y\rVert.
			\end{align*}
			where the last line uses the fact that $y_0,y$ are arbitrary. In particular, by setting $y_0:=y_n$ from step I, we get:
			$$\lVert \pst^*\rVert_{\infty,\xn}\leq |\pst^*(y_n)|+ \left(\max_{1\leq i\leq m}\lVert X_i\rVert\right)\sup_{y\in\xn}\lVert y_n-y\rVert\leq C_2\left(\max_{1\leq i\leq m}\lVert X_i\rVert\right),$$
			where $C_2:=3C_1$ with $C_1$ defined as specified in the proof of step I.
		\end{proof}
		The above lemma allows us to bound (with high probability) the $L^{\infty}$-norm of $\pst^*(\cdot)$ on $\xn$, using the tail assumption $\E\exp(t\lVert X_1\rVert^{\alpha})<\infty$ for some $t>0$ and $\alpha>0$. This is the focus of the next step. 
		
		\textbf{Step III}: For $K>0$, define the following two sets:
		$$A_{m,n,K}:=\left\{\int(\pst^*(u))^2\,d\nu(u)\geq K\right\},\quad \mathrm{and},$$
		$$\tilde{A}_{m,n,K}:=\left\{\lVert \pst^*\rVert_{\infty,\xn}\geq K\big(\log{n}\big)^{1/\alpha}\right\}.$$
		Then there exists $K_0>0$ such that for any $K\geq K_0$, we have:
		\begin{equation}\label{eq:nsmoothst3i}\lim_{m,n\to\infty}\P(\tilde{A}_{m,n,K})=0.\end{equation}
		and
		\begin{equation}\label{eq:nsmoothst3ii}\lim_{m,n\to\infty}\P(A_{m,n,K})=0.\end{equation}
		\begin{proof}[Proof of step III]
			By using the exponential Markov's inequality coupled with the standard union bound, we have:
			\begin{align*}
			\P\left(\max_{1\leq i\leq m} \lVert X_i\rVert\geq K(\log{m})^{1/\alpha}\right)&\leq m\P\left(\lVert X_1\rVert\geq K(\log{m})^{1/\alpha}\right)\\ &\leq m\exp(-tK^{\alpha}(\log{m}))\E\exp(t\lVert X_1\rVert^{\alpha})\overset{m\to\infty}{\longrightarrow}0
			\end{align*}
			provided $K>t^{-\alpha}$. Using the above observation coupled with step II,~\eqref{eq:nsmoothst3i} follows by choosing $K_0>C_2 t^{-\alpha}$.
			
			For the next part, we define another set:
			$$B_{m,n,\eps}:=\left\{\int \big|\pst^*(u)-\psm^*(u)|^2\,d\nu(u)\geq \eps\right\}$$
			for $\eps>0$, where, as in~\eqref{eq:dualf}, we have:
			$$W_2^2(\mu,\nu)=\int \lVert x\rVert^2\,d\mu(x)+\int \lVert y\rVert^2\,d\nu(y)-2\left(\int \psm(x)\,d\mu(x)+\int \psm^*(y)\,d\nu(y)\right).$$
			Now by using~\cite[Theorem 2.10]{del2019}, we have $\P(B_{m,n,\eps})\to 0$ as $m,n\to\infty$ for all $\eps>0$. As
			$$\int (\pst^*(u))^2\,d\nu(u)\leq 2\int \big|\pst^*(u)-\psm^*(u)|^2\,d\nu(u)+2\int (\psm^*(u))^2\,d\nu(u),$$
			\eqref{eq:nsmoothst3ii} follows with $K_0>2\int (\psm^*(u))^2\,d\nu(u)+1$ if we choose $\epsilon=1/2$.
		\end{proof}				
		
		We are now in a position to complete the proof of~\cref{thm:nsmooth} using steps I-III. Towards this direction, set $K':=2K_0$ where $K_0$ is defined as in the proof of step III  and observe that for any $M>0$,
		\begin{align}\label{eq:estimate1}
		&\;\;\;\limsup_{M\to\infty}\limsup_{m,n\to\infty}\P\left(\Bigg|\int \pst^*(u)\,d(\tnu-\nu)\Bigg|\geq Mr_d^{(n,n)}(\log{(1+m)})^{t_{d,\alpha}}\right)
		\nonumber \\&\leq \limsup_{M\to\infty}\limsup_{m,n\to\infty}\P\left(\Bigg|\int \pst^*(u)\,d(\tnu-\nu)\Bigg|\geq Mr_d^{(n,n)}(\log{(1+m)})^{t_{d,\alpha}}, A_{m,n,K'}^c\cap \tilde{A}_{m,n,K'}^c\right)\nonumber \\ &+\limsup\limits_{n\to\infty}\P(\tilde{A}_{m,n,K'})+\limsup\limits_{m,n\to\infty}\P(A_{m,n,K'})\nonumber \\&\leq \limsup_{M\to\infty}\limsup\limits_{m,n\to\infty}\P\left(\Bigg|\int \pst^*(u)\,d(\tnu-\nu)\Bigg|\geq Mr_d^{(n,n)}(\log{(1+m)})^{t_{d,\alpha}},A_{m,n,K'}^c\cap \tilde{A}_{m,n,K'}^c\right),
		\end{align}
		where the last step follows from step III. Observe that the left hand side of~\eqref{eq:estimate1} is the same as~\eqref{eq:nsmooth2}. Therefore, it is now enough to bound the right hand side of~\eqref{eq:estimate1}.
		
		In order to achieve the above task, let us define the following class of functions:
		$$\mathcal{C}^{\Gamma,L}(\xn):=\{f:\xn\to\R, f\ \mathrm{is}\ \mathrm{convex,} \ \lVert f\rVert_{\infty,\xn}\leq \Gamma, \ \lVert f\rVert_{L^2(\nu)}\leq L\}.$$
		By setting $\Gamma:=K'(\log{m})^{1/\alpha}$ and $L:=K'$,~\eqref{eq:estimate1} yields the following conclusion:
		\begin{align*}&\;\;\;\limsup_{M\to\infty}\limsup\limits_{m,n\to\infty}\P\left(\Bigg|\int \pst^*(u)\,d(\tnu-\nu)\Bigg|\geq Mr_d^{(n,n)}(\log{(1+m)})^{t_{d,\alpha}}\right)\\ & \leq \limsup_{M\to\infty}\limsup_{m,n\to\infty}\P\left(\sup_{f\in \mathcal{C}^{\Gamma,L}(\xn)}\Bigg|\int f\,d(\tnu-\nu)\Bigg|\geq Mr_d^{(n,n)}(\log{(1+m)})^{t_{d,\alpha}}\right).\end{align*}
		By an application of Markov's inequality, it thus suffices to show that:
		\begin{align}\label{eq:maintarget}
		\E\left[\sup_{f\in \mathcal{C}^{\Gamma,L}(\xn)}\Bigg|\int f\,d(\tnu-\nu)\Bigg|\right]=\mathcal{O}\left(r_d^{(n,n)}(\log{(1+m)})^{t_{d,\alpha}}\right).
		\end{align}
		In order to bound~\eqref{eq:maintarget}, we will use some standard empirical process techniques. In particular, by using~\cite[Theorem 5.11]{geer2000}, the following bound holds:
		\begin{align}\label{eq:dudleybd}
		&\;\;\;\E\left[\sup_{f\in \mathcal{C}^{\Gamma,L}(\xn)}\Bigg|\int f\,d(\tnu-\nu)\Bigg|\right]\nonumber \\ &\leq D\inf\left\{a\geq \frac{\Gamma}{\sqrt{n}}:a\geq \frac{D}{\sqrt{n}}\int_a^{\Gamma}\sqrt{\log N_{[]}(\eps,\mathcal{C}^{\Gamma,L}(\xn),L^2(\nu))}\,d\eps\right\},
		\end{align}
		for some positive constant $D>0$, where $N_{[]}(\eps,\mathcal{C}^{\Gamma,L}(\xn),L^2(\nu))$ is the $\eps$-bracketing number of the class of functions $\mathcal{C}^{\Gamma,L}(\xn)$ with respect to the $L^2(\nu)$ norm.
		Note that by~\cite[Equation 26]{kur2020convex}, we have:
		$$\log N_{[]}(\eps,\mathcal{C}^{\Gamma,L}(\xn),L^2(\nu))\leq \gamma_d\left(\log{\frac{\Gamma}{\eps}}\right)^{d+1}\left(\frac{L}{\eps}\right)^{d/2}$$
		for some $\gamma_d>0$ depending only on fand the diameter of $\xn$. 
		
		We will now bound the right hand side of~\eqref{eq:dudleybd}. Also we will use $D_d$ to denote changing constants which can depend on $d$.
		\begin{enumerate}
			\item \emph{When $d=1,2,3$}: Choose $a=D_d\frac{(\log{n})^{\frac{1}{\alpha}\vee \frac{2\alpha+2d\alpha-d+4}{4\alpha}}}{\sqrt{n}}$. Observe that:
			\begin{align*}
			\frac{1}{\sqrt{n}}\int_a^{\Gamma}\sqrt{\log N_{[]}(\eps,\mathcal{C}^{\Gamma,L}(\xn),L^2(\nu))}\,d\eps&\leq \frac{(\log{n})^{(d+1)/2}}{\sqrt{n}}\cdot\left[\frac{\eps^{1-d/4}}{1-d/4}\right]_{0}^{\Gamma}\\ &\lesssim \frac{(\log{n})^{(4-d)/(4\alpha)}\times (\log{n})^{(d+1)/2}}{\sqrt{n}}\lesssim a.
			\end{align*}
			\item \emph{When $d=4$}: Choose $a=D_d\frac{(\log{n})^{\frac{1}{\alpha}\vee \frac{7}{2}}}{\sqrt{n}}$. Observe that:
			\begin{align*}
			\frac{1}{\sqrt{n}}\int_a^{\Gamma}\sqrt{\log N_{[]}(\eps,\mathcal{C}^{\Gamma,L}(\xn),L^2(\nu))}\,d\eps&\leq \frac{(\log{n})^{5/2}}{\sqrt{n}}\cdot\left[\log{\eps}\right]_{D_d\Gamma/\sqrt{n}}^{\Gamma}\\ &\lesssim \frac{(\log{n})^{(7/2)}}{\sqrt{n}}\lesssim a.
			\end{align*}
			\item \emph{When $d>4$}: Choose $a=D_d\frac{(\log{n})^{2(1+d^{-1})}}{n^{2/d}}$. Observe that:
			\begin{align*}
			\frac{1}{\sqrt{n}}\int_a^{\Gamma_0}\sqrt{\log N_{[]}(\eps,\mathcal{C}^{\Gamma,L}(\xn),L^2(\nu))}\,d\eps&\leq \frac{(\log{n})^{(d+1)/2}}{\sqrt{n}}\cdot\left[\frac{\eps^{1-d/4}}{1-d/4}\right]_{a}^{\Gamma}\\ &\lesssim \frac{a^{1-d/4}(\log{n})^{(d+1)/2}}{\sqrt{n}}\lesssim a.
			\end{align*}
		\end{enumerate}
		This completes the proof after applying the same technique on the other $3$ terms on the right hand side of~\eqref{eq:nsmooth1}.
	\end{proof}
	
	\begin{proof}[Proof of~\cref{cor:fsam}]
		First observe that
		$$\E\left[\int \phmn^*\,d\hnu\right]=\E\left[\int \phmn^*\,d\bar{\nu}_m\right]=\int \phmn^*\,d\nu.$$
		Using the above observation and the same approach used as in the proof of~\cref{thm:nsmooth}, we will only focus on bounding
		\begin{equation}\label{eq:fsam2}
		\E\bigg|\int \pst^* \,d(\tnu-\nu)\bigg|.
		\end{equation}
		The general strategy to bound the term in~\eqref{eq:fsam2} is derived from some intermediate steps in the proofs of~\cite[Lemmas 3 and 4]{chizat2020faster}. We still present a sketch here for completeness.
		
		By the same argument as in the proof of~\cref{thm:nsmooth} and using the fact that there exists fixed $R>0$ such that $\max_{1\leq i\leq m}\lVert X_i\rVert\leq R$, we have $\pst^*(\cdot)$ is a convex and $R$-Lipschitz function on $\xn$. This observation implies:
		\begin{equation}\label{eq:fsam3}
		\E\bigg|\int \pst^* \,d(\tnu-\bar{\nu}_m)\bigg|\leq \E\left[\sup_{\psi\in\F_R(\xn)}\left|\int \psi \,d(\hnu-\nu)\right|\right]
		\end{equation}
		where $\F_R(\xn)$ is the set of convex and $R$-Lipschitz functions on $\xn$. By~\cite[Theorem 5.22]{Wainwright2019}, we then have:
		\begin{equation}\label{eq:fsam4}
		\E\left[\sup_{\psi\in\F_R(\xn)}\left|\int \psi \,d(\hnu-\nu)\right|\right]\lesssim \inf_{\delta>0}\left(\delta+n^{-1/2}\int_{\delta}^{R^2}\sqrt{\log{\cN_{\infty}(\F_R(\xn),\eps)}}\,d\eps\right),
		\end{equation}
		where $\cN_{\infty}(\F_R(\xn),\eps)$ is the $\eps$-covering number of the set $\F_R(\xn)$ with respect to the uniform metric. By using~\cite[Theorem 1]{guntuboyina2012} (also see~\cite{bronshtein1976varepsilon}), there exists constants $C_1,C_2>0$ such that whenever $\eps/R^2\leq C_1$, then  $\log{\cN_{\infty}(\F_R(\xn),\eps)}\leq C_2(u/R^2)^{-d/2}$. By using this bound in~\eqref{eq:fsam4}, we get:
		\begin{equation}\label{eq:fsam5}
		\E\left[\sup_{\psi\in\F_R(\xn)}\left|\int \psi \,d(\hnu-\nu)\right|\right]\lesssim \inf_{\delta>0}\left(\delta+n^{-1/2}\int_{\delta}^1 \eps^{-d/4}\,d\eps\right).
		\end{equation}
		Setting $\delta=0$ for $d<4$ and $\delta=n^{-2/d}$ for $d\geq 4$ in~\eqref{eq:fsam4}, followed by a direct application of~\eqref{eq:fsam3}, we have:
		$$\E\bigg|\int \pst^* \,d(\tnu-\bar{\nu}_m)\bigg|\lesssim r_d^{(n,n)}.$$
		This completes the proof.
	\end{proof}
	\begin{proof}[Proof of~\cref{thm:smoothwav}]
		For this proof, we will use an intermediate step in the proof of~\cref{thm:newubd}, which is~\eqref{eq:newubd7}, that can alternatively be written as:
		\begin{equation}\label{eq:smwav1}
		\E\left[\int \lVert \ttg(x)-\Tmn(x)\rVert^2\,d\tmu(x)\right]\lesssim\E|W_2^2(\tmu,\tnu)-W_2^2(\mu,\nu)|+\E\left|\int h(y)\,d(\tnu-\bar{\nu}_m)(y)\right|
		\end{equation}
		where $h(y):=\phmn^*(y)-(1/2)\lVert y\rVert^2$ and $C>0$ is some constant.
		As $\X$ and $\Y$ are compact sets, the function $h(\cdot)$ is Lipschitz. Therefore,
		$$\E\left|\int h(y)\,d(\tnu-\nu)(y)\right|\lesssim W_1(\tnu,\nu)\leq W_2(\tnu,\nu).$$
		Further, as $\Tmn(\cdot)$ is also Lipschitz, we further have:
		$$\E\left|\int h(y)\,d(\bar{\nu}_m-\nu)(y)\right|\lesssim W_1(\Tmn\#\tmu,\Tmn\#\mu)\lesssim W_1(\tmu,\mu)\leq W_2(\tmu,\mu).$$
		Finally, by the triangle inequality, we also have:
		\begin{align*}
		\E|W_2^2(\tmu,\tnu)-W_2^2(\mu,\nu)|&\lesssim \E|W_2(\tmu,\tnu)-W_2(\tmu,\nu)|+\E|W_2(\tmu,\nu)-W_2(\mu,\nu)|\\ &\leq \E W_2(\tmu,\mu)+\E W_2(\tnu,\nu).
		\end{align*}
		Combining the three displays above and plugging them back in~\eqref{eq:smwav1}, we get:
		$$\E\left[\int \lVert \ttg(x)-\Tmn(x)\rVert^2\,d\tmu(x)\right]\lesssim \E W_2(\tmu,\mu)+\E W_2(\tnu,\nu).$$
		The conclusion then follows from~\cite[Theorem 1]{weed2019estimation}.
	\end{proof}
	\begin{proof}[Proof of~\cref{thm:smooth}]
		\textbf{Part 1.} By the same arguments (see e.g.,~\eqref{eq:nsmooth1}) as used in the proof of~\cref{thm:nsmooth}, it suffices to show that
		\begin{equation}\label{eq:smooth1}
		\E\left|\int \pst^*(u)(\tfn^{M'}(u)-\fn(u)) \,du\right|\lesssim r_{d,s}^{(n,n)}
		\end{equation}
		for some $M'>0$.
		
		The general structure of the proof is similar to that of~\cref{thm:nsmooth}. The crucial observation is that $\tfm^{M'}(\cdot)$ and $\tfn^{M'}(\cdot)$ are elements of $C^s(\X;TM)$ and $C^s(\Y;TM)$ respectively, for any $M'>0$. Note that, by Caffarelli regularity theory; see~\cite[Theorem 33]{hutter2021minimax}, there exists $M'>0$ such that $\lVert \pst^*(\cdot)\rVert_{C^{s+2}(\Y)}\leq M'.$ 
		
		Next, let us define the following class of functions:
		$$\G_t^L(\Y):=\{g:\xn\to\R,\ g(\cdot)\ \mbox{is}\ \mbox{convex},\ \lVert g\lVert_{C^t(\xn)}\leq L\}.$$
		
		Observe that				
		\begin{align}\label{eq:smoothnew1}
		\E\left|\int \pst^*(u)(\tfn^{M'}(u)-\fn(u)) \,du\right|&\leq \E\sup_{g\in \G_{s+2}^{M'}(\Y)}\left|\int g(u)(\tfn^{M'}(u)-\fn(u)) \,du\right|\nonumber \\ &\leq 2\E\sup_{g\in \G_{s+2}^{M'}(\Y)}\left|\int g(u)(\hfn(u)-\fn(u)) \,du\right|+r_{d,s}^{(n,n)}
		\end{align}	
		where the last line follows from~\eqref{eq:kerproj}.
		
		Set $K_{d,h_n}(\cdot):=h_n^{-d}K_d(\cdot/h_n)$. Following the same decomposition as in~\cite{Radulovic2000}, we write:
		\begin{align}\label{eq:smooth4}
		&\;\;\;\;\E\sup_{g\in \G_{s+2}^{M'}(\xn)}\left|\int g(u)(\hfn(u)-\fn(u)) \,du\right|\nonumber \\ &=\E\sup_{g\in \G_{s+2}^{M'}(\xn)}\left|\int g(u+u')K_{d,h_n}(u')\,d\hnu(u)\,du'-\int g(u)\fn(u)\,du\right|\nonumber \\ &\leq \E\sup_{g\in \G_{s+2}^{M'}(\xn)}\left|\int g(u+u')K_{d,h_n}(u')\,d(\hnu-\nu)(u)\,du'\right|\nonumber \\ &\;\;\;\;\;\;\;\;\;\;+\sup_{g\in \G_{s+2}^{M'}(\xn)}\left|\int g(u+u')K_{d,h_n}(u')\fn(u)\,du\,du'-\int g(u)\fn(u)\,du\right|.
		\end{align}
		We will now bound the two terms on the right hand side of~\eqref{eq:smooth4}. For the first term, define 
		$$\bar{g}_n(u):=\int g(u+u')K_{d,h_n}(u)\,du'.$$				
		
		If $g\in \G_t^L(\xn^o)$, then by~\cite[Proposition 8.10]{folland1999real} and using~\cref{as:kernel}, we have $\bar{g}_n\in \G_{s+2}^{cM'}(\xn^o)$, for some constant $c>0$ (depending on the constants involved in~\cref{as:kernel} and the diameter of $\xn$). Combining these observations with~\eqref{eq:smooth4}, we get:
		\begin{align}\label{eq:smoothnew2}
		&\;\;\;\;\E\sup_{g\in \G_{s+2}^{M'}(\xn)}\left|\int g(u+u')K_{d,h_n}(u')\,d(\hnu-\nu)(u)\,du'\right|\nonumber \\ &\leq \E\sup_{g\in \G_{s+2}^{cM'}(\xn)}\left|\int g(u)\,d(\hnu-\nu)(u)\right|\nonumber \\ &\leq D\inf\left\{a\geq \frac{cM'}{\sqrt{n}}:a\geq \frac{D}{\sqrt{n}}\int_a^{cM'}\sqrt{\log N_{[]}(\eps,\G_{s+2}^{cM'}(\xn),L^2(\nu))}\,d\eps\right\},
		\end{align}
		for some positive constant $D>0$, where $N_{[]}(\eps,\G_{s+2}^{cM'}(\xn),L^2(\nu))$ is the $\eps$-bracketing entropy of the class of functions $\G_{s+2}^{cM'}(\xn)$ with respect to the $L^2(\nu)$ norm. The last line follows from standard empirical process theory as used in the proof of~\cref{thm:nsmooth}; see~\eqref{eq:dudleybd}. 
		Note that by~\cite[Corollary 2.7.2]{vaart2013}, we have:
		$$\log N_{[]}(\eps,\G_{s+2}^{cL}(\xn^o),L^2(\nu))\leq \gamma_d\left(\frac{1}{\eps}\right)^{d/(s+2)}$$
		for some $\gamma_d>0$ depending only on dimension and the diameter of $\xn$. 
		
		We now plug-in the above bound into~\eqref{eq:smoothnew2}. By using $D_d$ to denote constants that change with $d$ and choosing $a=D_d n^{-1/2}$ for $2(s+2)>d$, $a=D_d n^{-1/2}\log{(1+n)}$ for $2(s+2)=d$ and $a=D_d n^{-(s+2)/d}$ for $2(s+2)<d$ in~\eqref{eq:smoothnew2}, we have:
		\begin{equation}\label{eq:smoothnew3}
		\E\sup_{g\in \G_{s+2}^{M'}(\xn)}\left|\int g(u+u')K_{d,h_n}(u')\,d(\hnu-\nu)(u)\,du'\right|\lesssim r_{d,s}^{(n,n)}.
		\end{equation} 
		
		We now move on to the bounding the second term on the right hand side of~\eqref{eq:smooth4}. For this part, our main technical tool will be the classical arguments for smoothed empirical processes developed in~\cite{Gine2008}. Towards this direction, set $\bar{g}(u)=g(-u)$ (different from $\bar{g}_n(\cdot)$ defined earlier) for $g(\cdot)\in \G_{s+2}^{M'}$ and note that by~\cite[Lemma 4]{Gine2008}, we have:
		\begin{align}\label{eq:smoothnew4}
		\left|\int g(u+u')K_{d,h_n}(u')\fn(u)\,du\,du'-\int g(u)\fn(u)\,du\right|=\left|\int K_d(u)\left[(\bar{g}*\fn)(h_n u)-(\bar{g}*\fn)(0)\right]\,du\right|,
		\end{align}
		where $(\bar{g}*\fn)(\cdot)$ is the standard convolution between $\bar{g}(\cdot)$ and $\fn(\cdot)$, and with a notational abuse $0$ denotes the $d$-dimensional zero vector. The important observation now is to note that $(\bar{g}*\fn)(\cdot)$ belongs to a higher order Sobolev class compared to $\bar{g}(\cdot)$ and $\fn(\cdot)$. In particular, as $\fn(\cdot)\in C^s(\Y;M)$ and $\bar{g}(\cdot)\in \G_{s+2}^{M'}(\Y)$, we have $(\bar{g}*\fn)(\cdot)\in \G_{2s+2}^{M''}(\Y)$ where $M''$ depends on both $M'$ and $M$.
		
		Next, write $D^t(\bar{g}*\fn)(\cdot)$ to be the $t$-th derivative of $(\bar{g}*\fn)(\cdot)$ and note that by a multivariate Taylor's approximation
		\begin{align*}
		&\;\;\;\;\int K_d(u)\left[(\bar{g}*\fn)(h_n u)-(\bar{g}*\fn)(0)\right]\,du\\ &=\int K_d(u)\sum_{r=1}^{2s+1} h_n^r\sum_{(i_1,i_2,\ldots ,i_r)\in \{1,2,\ldots ,d\}^r} [D^r(\bar{g}*\fn)(0)]_{i_1,\ldots ,i_r} u_{i_1}\ldots u_{i_r}\,du+O(h_n^{2s+2}).
		\end{align*}
		Recall that $K_d(u)=K(u_1)K(u_2)\ldots K(u_d)$. As $K(\cdot)$ is of order $2s+2$ (see~\cref{as:kernel}), all the integrals on the right hand side of the above display vanish. We then appeal to~\eqref{eq:smoothnew4} to get:
		$$\left|\int g(u+u')K_{d,h_n}(u')\fn(u)\,du\,du'-\int g(u)\fn(u)\,du\right|\lesssim h_n^{2s+2}=\lesssim n^{-\frac{2s+2}{d+2s}}(\log{n})^{2s+2}.$$
		We now compare the right hand side of the above display with $r_{d,s}^{(n,n)}$.
		
		\textbf{When $d<2(s+2)$:} $d+2s<4(s+1)$, and therefore $\frac{2s+2}{d+2s}>\frac{1}{2}$. This implies $n^{-\frac{2s+2}{d+2s}}(\log{n})^{2s+2}\lesssim n^{-\frac{1}{2}}=r_{d,s}^{(n,n)}$.
		
		\textbf{When $d=2(s+2)$:} In this case $n^{-\frac{2s+2}{d+2s}}(\log{n})^{2s+2}\lesssim n^{-\frac{1}{2}}(\log{n})^{2s+2}\lesssim r_{d,s}^{(n,n)}$.
		
		\textbf{When $d>2(s+2)$:} Note that
		$$\frac{2s+2}{d+2s}>\frac{s+2}{d}\Leftrightarrow 2ds+2d>sd+2s^2+2d+4s\Leftrightarrow d>2(s+2).$$
		Therefore, once again $n^{-\frac{2s+2}{d+2s}}(\log{n})^{2s+2}\lesssim n^{-\frac{s+2}{d}}=r_{d,s}^{(n,n)}$.
		
		Therefore, combining the above observations, we have:
		$$\left|\int g(u+u')K_{d,h_n}(u')\fn(u)\,du\,du'-\int g(u)\fn(u)\,du\right|\lesssim r_{d,s}^{(n,n)}.$$
		Combining the above display with~\eqref{eq:smoothnew3} establishes~\eqref{eq:smooth1}. 
		
		\textbf{Part 2.} This proof uses ideas from~\cite[Theorem 7]{Hansen08},~\cite{Parzen1962}~and~\cite[Lemmas 2 and 3]{arias2016estimation}. First recall all the notation introduced in~\cref{def:holder}. Next, we will prove the following sequence of displays:
		\begin{equation}\label{eq:s1claim1}
		\limsup\limits_{m,n\to\infty}\max_{k\leq s}\max_{|\mm|=k}\lVert\partial^{\mm}\E\hfm\rVert_{L^{\infty}(\tilde{\X})}\leq (T-1)M,
		\end{equation}
		\begin{equation}\label{eq:s1claim1p5}
		\limsup_{m,n\to\infty}\ \lVert \E\hfm-\fm\rVert_{L^{\infty}(\tilde{\X})}=0
		\end{equation}
		\begin{equation}\label{eq:s1claim2}
		\limsup_{m,n\to\infty}\P\left(\lVert \hfm-\E\hfm\rVert_{C^{s}(\tilde{\X})}\geq \eps\right)=0,
		\end{equation}
		for any arbitrary $\eps>0$ and $L^{\infty}(\X)$ denotes the uniform norm on $\X$.
		
		Clearly,~\eqref{eq:s1claim1},~\eqref{eq:s1claim1p5},~and~\eqref{eq:s1claim2} together yield part 1 of the theorem.
		
		\emph{Proof of~\eqref{eq:s1claim1}.} Observe that
		\begin{equation}\label{eq:s11}
		\E\hfm(x)=\frac{1}{h_m^d}\E K_d\left(\frac{x-X_1}{h_m}\right)=\frac{1}{h_m^d}\int K_d\left(\frac{x-z}{h_m}\right)\fm(z)\,dz.
		\end{equation}	
		Since the maximums taken in~\eqref{eq:s1claim1} are over finite sets, it suffices to show that for any fixed $\mm$ with $|\mm|\leq s$, we have:
		\begin{equation}\label{eq:s12}
		\sup_{x\in\tilde{\X}}|\partial^{\mm}\E\hfm(x)|=\sup_{x\in\X}\bigg|\frac{1}{h_n^d}\int K_d\left(\frac{z}{h_n}\right)\partial^{\mm} f(x+z)\,dz\Bigg|\leq (T-1)M.
		\end{equation}
		Here the first equality in the above display follows from~\eqref{eq:s11} and Fubini's Theorem. Here $\partial^{\mm}\fm(\cdot)$ is defined in the weak sense, i.e., it is defined naturally in the interior of the support of $\fm(\cdot)$, denoted by $\X$; it is set to be $0$ outside $\X$ and defined arbitrarily on the boundary of $\X$. Note that the definition on the boundary doesn't matter as we are integrating with respect to the Lebesgue measure and the boundary of $\X$ has Lebesgue measure $0$.
		
		Next note that, by~\eqref{eq:s12}, we have:
		$$\sup_{x\in\tilde{\X}}|\partial^{\mm}\E\hfm(x)|\leq \lVert \fm\rVert_{C^s(\X)}h_m^{-d}\int |K_d(z/h_m)|\,dz\leq (T-1)\lVert\fm\rVert_{C^s(\X)}.$$
		This establishes~\eqref{eq:s1claim1}.				
		
		\emph{Proof of~\eqref{eq:s1claim1p5}.} First note that, as $\tilde{\X}$ is a compact subset of $\X^o$, there exists $\delta>0$ such that
		$$\tilde{\X}_{\delta'}:=\{x+z:\ \lVert z\rVert\leq \delta;,\ x\in\tilde{\X}\}\subseteq \X^o\qquad \forall \ 0<\delta'\leq \delta.$$ 
		Clearly, $\tilde{\X}_{\delta'}$ is compact for all $\delta'>0$. Fix an arbitrary $\delta'\leq \delta$.
		By using~\eqref{eq:s11} and a change of variable formula, we have:
		\begin{align*}
		\lVert \E\hfm-\fm\rVert_{L^{\infty}(\tilde{\X})}&=\sup_{x\in\tilde{\X}}\bigg|\frac{1}{h_m^d}\int K_d\left(\frac{z}{h_m}\right)(f(x+z)-f(x))\,dz\bigg| \\& \leq (T-1)\sup_{x\in\tilde{\X}}\sup_{\lVert z\rVert\leq \delta'}|f(x+z)-f(x)|+2M\int_{\lVert z\rVert>\delta' h_m^{-1}} |K_d(z)|\,dz\\ & \leq (T-1)M\delta'+2M\left(\frac{h_m}{\delta'}\right)^{2s+2}\int \lVert z\rVert^{2s+2}|K_d(z)|\,dz.
		\end{align*}
		Observe that as $m,n\to\infty$, the second term on the right hand side of the above display converges to $0$. This implies
		$$\limsup\limits_{m,n\to\infty}\lVert \E\hfm-\fm\rVert_{L^{\infty}(\tilde{\X})}\leq (T-1)M\delta'.$$
		As $\delta'$ can be chosen arbitrarily small, this completes the proof of~\eqref{eq:s1claim1p5}.
		
		\emph{Proof of~\eqref{eq:s1claim2}.} The main technical tool for this part is~\cref{lem:emproc} which we borrow from~\cite[Lemma 9]{arias2016estimation} (also see~\cite[Theorem 4.1]{Mason2012}). The proof is very similar to~\cite[Lemma 3]{arias2016estimation}.  Consider the following class of functions:
		$$\mathcal{G}=\left\{g_x(z,h):\ g_x(z,h)=\partial^{\mm}K_d\left(\frac{(x-z)}{h}\right),\ x\in\tilde{\X},\ |\mm|\leq s \right\}.$$
		Observe that
		$$\sup_{|\mm|\leq s}\sup_{x\in\tilde{\X}}\sup_{h\in (0,1)} h^{-d}  \E\left[\partial^{\mm}K_d\left(\frac{(x-z)}{h}\right)\right]^2\leq \lVert K\rVert_{C^s(\R^d)}\lVert f\rVert_{C^s(\X)}\sup_{|\mm|\leq s}\int |\partial^{\mm} K_d(v)|\,dv<\infty.$$
		Further, by~\cref{as:kernel}, $\partial^{\mm} K_d(\cdot)$ is differentiable for each $|\mm|\leq s$. Consequently $\mathcal{G}$ is point wise measurable and of VC-type (see~\cite[Lemma A.1]{Einmahl2000}; also see~\cite[Section 2.6]{vaart2013} for definitions of point wise differentiability and VC classes). This verifies the assumptions of~\cref{lem:emproc}. Observe that
		$$\frac{1}{n}\sum_{i=1}^m \partial^{\mm} K\left(\frac{x-X_i}{h_m}\right)=h_m^{d+|\mm|}\partial^{\mm} \hfm(x),\qquad \E\left[\partial^{\mm}K\left(\frac{x-X}{h_m}\right)\right]=h_m^{d+|\mm|}\E\left[\partial^{\mm}\hfm(x)\right].$$
		A direct application of~\cref{lem:emproc} for all $|\mm|\leq s$,  then implies
		$$\sup_{x\in\tilde{\X}}\sqrt{\frac{m}{h_m^d\log{m}}}\cdot h_m^{s+d} \lVert \hfm-\E\hfm\rVert_{C^s(\tilde{\X})}=O_p(1).$$
		Using the observation that $mh_m^{d+2s}/\log{m}\to 0$ as $m\to\infty$ then completes the proof.				
	\end{proof}
	\begin{proof}[Proof of~\cref{prop:wassrate}]
		As $\mu\neq \nu$, we have $W_2(\mu,\nu)>0$. Therefore,
		$$|W_2(\tmu,\tnu)-W_2(\mu,\nu)|=\frac{W_2^2(\tmu,\tnu)-W_2^2(\mu,\nu)|}{W_2(\tmu,\tnu)+W_2(\mu,\nu)}\leq \frac{W_2^2(\tmu,\tnu)-W_2^2(\mu,\nu)|}{W_2(\mu,\nu)}.$$
		The conclusion then follows from~\cref{thm:smooth}.
	\end{proof}
	
	\begin{proof}[Proof of~\cref{thm:dismooth}]
		Recall that $\tmu$ and $\tnu$ are defined as the empirical distributions induced by $M=n^{\frac{s+2}{2}}$ random samples drawn from $\hfm$ and $\hfn$ respectively, where $\hfm$, $\hfn$ are the kernel density estimates as presented in~\eqref{eq:den}. Let us write $\mu_{h_n}$ and $\nu_{h_n}$ for the probability measure induced by the kernel density estimates $\hfm$ and $\hfn$ respectively. Once again, by using~\cref{thm:nsmooth},~\eqref{eq:newubd8}, it suffices to prove the following:
		\begin{equation}\label{eq:dismoothnew1}
		\E\left|W_2^2(\tmu,\tnu)-W_2^2(\mu,\nu)\right|.
		\end{equation}
		Next note that by the triangle inequality,~\eqref{eq:dismoothnew1} can be bounded above by:
		\begin{align}\label{eq:dismoothnew2}
		\E\big|W_2^2(\tmu,\nu_{h_n})&-W_2^2(\tmu,\tnu)\big|+\E\big|W_2^2(\tmu,\nu_{h_n})-W_2^2(\mu_{h_n},\nu_{h_n})\big|\nonumber \\&+\E\big|W_2^2(\mu_{h_n},\nu_{h_n})-W_2^2(\tmu,\tnu)\big|.
		\end{align}
		Next note that, by~\cref{thm:smooth}, we have:
		\begin{equation}\label{eq:dismoothnew3}
		\E\left|W_2^2(\mu_{h_n},\nu_{h_n})-W_2^2(\tmu,\tnu)\right|\lesssim r_{d,s}^{(n,n)}.
		\end{equation}
		Next we show that				
		\begin{equation}\label{eq:dismooth1}
		\E\left|W_2^2(\tmu,\nu_{h_n})-W_2^2(\mu_{h_n},\nu_{h_n})\right|\lesssim r_{d,s}^{(n,n)}.
		\end{equation}				
		The other term in~\eqref{eq:dismoothnew2} can be bounded similarly.
		
		Note that, conditioned on $X_1,\ldots ,X_n,Y_1,\ldots ,Y_n$, $\mu_{h_n}$ and $\nu_{h_n}$ are non-random measures and $\tmu$ and $\tnu$ are the empirical distributions on $M=n^{\frac{s+2}{2}}$ random samples from the measures $\mu_{h_n}$ and $\nu_{h_n}$, respectively. Therefore, conditioned on $X_1,\ldots ,X_n,Y_1,\ldots ,Y_n$ (which have fixed compact supports), we can invoke~\cref{cor:fsam} to get:
		$$\E\left|W_2^2(\tmu,\nu_{h_n})-W_2^2(\mu_{h_n},\nu_{h_n})\right| \lesssim r_{d}^{(M,M)},$$
		with $M=n^{\frac{s+2}{2}}$. Recall that:
		$$r_d^{(M,M)}=\begin{cases} n^{-\frac{s+2}{4}} & \mbox{if}\ d\leq 3\\ n^{-\frac{s+2}{4}}\log{(1+n)} & \mbox{if}\ d=4\\ n^{-\frac{s+2}{d}} & \mbox{if}\ d>4\end{cases}.$$
		It therefore only remains to compare $r_d^{(M,M)}$ and $r_{d,s}^{(n,n)}$.
		
		\textbf{Case 1:} $d\leq 2(s+2)$. In this case, if $d=1,2,3$, then $r_d^{(M,M)}=n^{-\frac{s+2}{4}}=n^{-\frac{1}{2}}\times n^{-\frac{s}{4}}\lesssim n^{-\frac{1}{2}}$. If $d=4$, then $r_d^{(M,M)}= n^{-\frac{s+2}{4}}\log{(1+n)}=n^{-\frac{1}{2}}\times \left(n^{-\frac{s}{4}}\log{n}\right)\lesssim n^{-\frac{1}{2}}$. If $d>4$, then $r_d^{(M,M)}=n^{-\frac{s+2}{d}}\lesssim n^{-\frac{1}{2}}$ as $\frac{s+2}{d}\geq \frac{1}{2}$. Therefore, in all the cses, $r_d^{(M,M)}\lesssim n^{-\frac{1}{2}}=r_{d,s}^{(n,n)}$ for $d\leq 2(s+2)$.
		
		\textbf{Case 2:} $d>2(s+2)$. As $s>0$, then $d>4$. In this case, once again $r_d^{(M,M)}=n^{-\frac{s+2}{d}}=r_{d,s}^{(n,n)}$.
		
		This establishes~\eqref{eq:dismooth1} and completes the proof.
	\end{proof}		
	\subsection{Proofs from~\cref{sec:App}}\label{sec:Appf}
	\begin{proof}[Proof of~\cref{thm:barate}]
		First define the following measure:
		\begin{equation*}
		\bmor:=\left(\frac{1}{2}\mbox{Id}+\frac{1}{2}\Tmn\right)\#\tmu.
		\end{equation*}
		Fix any $\gamma\in\tgm$. By applying the triangle inequality followed by a power mean inequality, we have:
		\begin{equation}\label{eq:barate1}
		\sup_{\gamma\in\tgm}W_2^2\big(\Bmnz^{\gamma},\Bmn\big)\lesssim W_2^2\big(\bmor,\Bmn\big)+\sup_{\gamma\in\tgm} W_2^2\big(\Bmnz^{\gamma},\bmor\big).
		\end{equation}
		Next observe that $\bmor$ is the empirical distribution corresponding to $m$ random samples drawn according to $\Bmn$. Therefore, by using~\cite[Theorem 1]{Fournier2015}, we get:
		\begin{equation}\label{eq:barate2} W_2^2\big(\bmor,\Bmn\big)\lesssim r_d^{(m,m)}.\end{equation}
		Next we will bound the second term on the right hand side of~\eqref{eq:barate1}. Towards this direction, recall the definition of $\Pi(\cdot,\cdot)$ from~\cref{sec:setting}. Consider the following coupling:
		$$\pi_0^{\gamma}:=\left(\frac{1}{2}\mbox{Id}+\frac{1}{2}\ttg,\frac{1}{2}\mbox{Id}+\frac{1}{2}\Tmn\right)\#\tmu.$$
		Observe that $\pi_0^{\gamma}\in \Pi\big(\Bmnz^{\gamma},\bmor\big)$.  By plugging the coupling $\pi_0^{\gamma}$ into the definition of $2$-Wasserstein distance in~\eqref{eq:wass}, we further get:
		\begin{align}\label{eq:barate3}
		\sup_{\gamma\in\tgm} W_2^2\big(\Bmnz^{\gamma},\bmor\big)&\leq \sup_{\gamma\in\tgm} \int \lVert x-y\rVert^2\,d\pi_0^{\gamma}(x,y)\nonumber \\ &=\sup_{\gamma\in\tgm} \int \lVert \ttg(x)-\Tmn(x)\rVert^2\,d\tmu(x)\nonumber \\ &= O_p\left(r^{(m,n)}_d\times (\log{(1+\max\{m,n\})})^{t_{d,\alpha}}\right)
		\end{align}
		where the last inequality follows from~\cref{thm:nsmooth}.  Combining~\eqref{eq:barate2}~and~\eqref{eq:barate3}  with~\eqref{eq:barate1} completes the proof.
	\end{proof}
	\begin{proof}[Proof of~\cref{thm:dethresh}]
		Let $\oto(\cdot)$ and $\ott(\cdot)$ be the optimal transport maps from $\mu^{(n)}$ to $\upsilon_1$ and $\nu^{(n)}$ to $\upsilon_2$. Set  $$\hxr_{ij}:=K_1(\oto(X_i),\oto(X_j)),\qquad   \hyr_{ij}:=K_2(\ott(Y_i),\ott(Y_j))$$ and define the oracle version of $\hhs$ as follows: 
		\begin{equation}\label{eq:dethresh1}\ohs:=\underbrace{n^{-2}\sum_{i,j}\hxr_{ij}\hyr_{ij}}_{\haor}+\underbrace{n^{-4}\sum_{i,j,r,s}\hxr_{ij}\hyr_{rs}}_{\hadr}-2\underbrace{n^{-3}\sum_{i,j,r} \hxr_{ij}\hyr_{ir}}_{\harr}.\end{equation}
		The proof of~\cref{thm:dethresh} now proceeds using the following steps:
		
		\textbf{Step I:} We show that:
		\begin{equation}\label{eq:dethresh2}
		\E\big|\ohs-\nhs(\pi^{(n)}|\mu^{(n)}\times \nu^{(n)})\big|\lesssim n^{-1/2},
		\end{equation}
		where $\ohs(\cdot|\cdot)$ is defined in~\eqref{eq:popmeas}.
		
		\textbf{Step II:} We prove that:
		\begin{equation}\label{eq:dethresh3}	
		\E\big|\ohs-\hhs\big|\lesssim \sqrt{r_d^{(n,n)}}.
		\end{equation}
		\textbf{Step III:} We combine steps I and II to prove~\cref{thm:dethresh}. Let us begin with this step first. Note that by using the triangle inequality, we have:
		\begin{equation}\label{eq:dethresh4}
		\hhs\geq \nhs(\pi^{(n)}|\mu^{(n)}\times \nu^{(n)})-\big|\ohs-\nhs\big|-\big|\ohs-\hhs\big|.
		\end{equation}
		Next observe that by steps I and II,
		$$\max\bigg\{\big|\ohs-\nhs\big|,\big|\ohs-\hhs\big|\bigg\}=O_p\big(\sqrt{r_{d_1,d_2}^{(n,n)}}\big).$$
		Using the above display with~\eqref{eq:dethresh4} and the assumption $(r_{d_1,d_2}^{(n,n)})^{-1/2}\nhs(\pi^{(n)}|\mu^{(n)}\times \nu^{(n)})\to\infty$, we have:
		$$\big(r_{d_1,d_2}^{(n,n)}\big)^{-1/2}\hhs\overset{P}{\longrightarrow}\infty.$$
		Therefore, as $n\sqrt{r_{d_1,d_2}^{(n,n)}}\to\infty$ and $c_{n,\alpha}=O(1)$ (see~\cite[Theorem 4.1]{deb2021multivariate}), we have:
		$$\E\phi_{n,\alpha}=\P(n\times\hhs\geq c_{n,\alpha})\to 1$$
		under $(r_{d_1,d_2}^{(n,n)})^{-1/2}\nhs(\pi^{(n)}|\mu^{(n)}\times \nu^{(n)})\to\infty$. This completes the proof.
		
		It therefore remains to prove steps I and II. For step I, let $(X_1',Y_1'),\ldots ,(X_n',Y_n')\overset{i.i.d.}{\sim}\pi^{(n)}$. Fix an arbitrary $1\leq j\leq n$. Let $\hat{A}_{n,1,j}^{\mathrm{OR},'}$ be the same as $\haor$ except with $(X_j,Y_j)$ replaced by $(X_j',Y_j')$. It is easy to check by the compactness of supports of all distributions involved, that:
		$$\max\limits_{1\leq j\leq n}\big|\haor-\hat{A}_{n,1,j}^{\mathrm{OR},'}\big|\lesssim n^{-1}.$$
		Therefore by using Mcdiarmid's inequality (see~\cite[Theorem 6.5]{Boucheron2013}), we have, for any $t>0$, 
		$$\P\left(\sqrt{n}(\haor-\E\haor)\geq t\right)\leq \exp(-Ct^2)$$
		for some constant $C>0$ free of $n$ and $t$. Similar concentrations can be derived for $\hadr$ and $\harr$. Combining these concentrations with the observation that
		$$\nhs(\pi^{(n)}|\mu^{(n)}\times \nu^{(n)})=\E\haor+\E\hadr-2\E\harr$$
		completes the proof of step I. 
		
		We now move on to step II. Recall the definition of $\hhs$ from~\eqref{eq:emprhsic} and write:
		$$\hhs=\underbrace{n^{-2}\sum_{i,j}\hx_{ij}\hy_{ij}}_{\hao}+\underbrace{n^{-4}\sum_{i,j,r,s}\hx_{ij}\hy_{rs}}_{\had}-2\underbrace{n^{-3}\sum_{i,j,r} \hx_{ij}\hy_{ir}}_{\har}.$$
		By the Lipschitzness of $K_1(\cdot,\cdot)$ and $K_2(\cdot,\cdot)$, we have:
		$$\big|\hx_{ij}-\hxr_{ij}\big|\lesssim \lVert \hto(X_i)-\oto(X_i)\rVert+\lVert \hto(X_j)-\oto(X_j)\rVert,$$
		$$\big|\hy_{ij}-\hyr_{ij}\big|\lesssim \lVert \htt(Y_i)-\ott(Y_i)\rVert+\lVert \htt(Y_j)-\ott(Y_j)\rVert.$$
		Therefore, by using the fact that the probability measures $\upsilon_1$ and $\upsilon_2$ are compactly supported, we get:
		$$\big|\hao-\haor\big|\lesssim \frac{1}{n}\sum_{i=1}^n \lVert \hto(X_i)-\oto(X_i)\rVert+\frac{1}{n}\sum_{j=1}^n \lVert \htt(Y_j)-\ott(Y_j)\rVert.$$
		The same bound can similarly be verified for $|\had-\hadr|$ and $|\har-\harr|$. Combining these observations, we have:
		\begin{align*}
		\big|\hhs-\ohs\big|&\lesssim \frac{1}{n}\sum_{i=1}^n \lVert \hto(X_i)-\oto(X_i)\rVert+\frac{1}{n}\sum_{j=1}^n \lVert \htt(Y_j)-\ott(Y_j)\rVert\nonumber \\ &\leq \sqrt{\frac{1}{n}\sum_{i=1}^n \lVert \hto(X_i)-\oto(X_i)\rVert^2}+\sqrt{\frac{1}{n}\sum_{j=1}^n \lVert \htt(Y_j)-\ott(Y_j)\rVert^2}.
		\end{align*}
		Step II then follows by invoking~\cref{cor:fsam}.
	\end{proof}
	\section{Auxiliary definitions and results}\label{sec:mdefpt}
	\begin{defn}[Subdifferential set and subgradient]\label{def:subg}
		Given a convex function $f:\R^d\to \R\cup\{\infty\}$, we define the \emph{subdifferential} set of $f(\cdot)$ at $x\in\mbox{dom}(f):=\{z\in\R^d:f(z)<\infty\}$ as follows:
		$$\partial f(x):=\{\xi\in\R^d:\ f(x)+\langle \xi, y-x\rangle\leq f(y),\quad \mbox{for}\ \mbox{all}\ y\in\R^d\}.$$
		Any element in the set $\partial f(x)$ is called a \emph{subgradient} of $f(\cdot)$ at $x$.
	\end{defn}
	\begin{defn}[Strong convexity]\label{def:strongcon}
		A function $f:\R^d\to\R\cup\{\infty\}$ is strongly convex with parameter $\lambda>0$, if, for all $x,y\in\mbox{dom}(f)=\{z\in\R^d:f(z)<\infty\}$, the following holds:
		$$f(y)\geq f(x) + \langle \xi_x,y-x\rangle+\frac{\lambda}{2}\lVert y-x\rVert^2,$$
		where $\xi_x\in \partial f(x)$, the subgradient of $f(\cdot)$ at $x$ as in~\cref{def:subg}.
	\end{defn}						
		
	\begin{defn}[Wavelet basis]\label{def:wavelet}
		We present our main assumptions on the wavelet basis discussed in~\cref{sec:wavelet} only for the wavelets on the space $\X$. The same assumptions are also required for the wavelets on $\Y$. These are essentially a subset of the assumptions laid out in~\cite[Appendix E]{weed2019estimation} as we heavily rely on~\cite[Theorem 1]{weed2019estimation} for proving~\cref{thm:smoothwav}.
		
		\begin{enumerate}
			\item \textbf{(Regularity).} Fix $r>\max\{s,1\}$. The functions in $\bph$ and $\bps_j$, $j\geq 0$ have $r$ continuous derivatives, and all polynomials of degree at most $r$ on $\X$ lie in the span of the functions in $\bph$.
			\item \textbf{(Tensor construction).} Each $\psi(\cdot)\in\bps_j$ can be expressed as $\psi(x)=\prod_{i=1}^d \psi_i(x_i)$, where $x=(x_1,\ldots ,x_d)$, for some univariate functions $\psi_i(\cdot)$'s.
			\item \textbf{(Locality).} For each $\psi(\cdot)\in\bps_j$ there exists a rectangle $I_{\psi}\subseteq \X$ such that $\mbox{supp}(\psi)\subseteq I_{\psi}$, $\mbox{diam}(I_{\psi})\leq C_1\cdot 2^{-j}$, and $\sup_{x\in\X}\sum_{\psi(\cdot)\in\bps_j}\mathbbm{1}(x\in I_{\psi})\leq C_2$ for some constants $C_1,C_2>0$.
			\item \textbf{(Bernstein estimate).} $\lVert \nabla f\rVert_{L^2(\X)}\leq C_3\cdot 2^j\lVert f\rVert_{L^2(\X)}$ for any $f(\cdot)$ in the span of the functions in $\mbox{span}\left(\bph\cup\left\{\cup_{0\leq k<j}\bps_j\right\}\right)$. Here $C_3$ is some positive constant.
		\end{enumerate}
	\end{defn}
			
	\begin{lemma}[Strong convexity and Lipschitzness, see~\cite{hiriart1993}]\label{lem:assm}
		$\phmn^*(\cdot)$ is strongly convex with parameter $(1/L)$ if and only if $\Tmn(\cdot)$ is $L$-Lipschitz continuous.
	\end{lemma}
	\begin{lemma}[Gradient of dual]\label{lem:gradual}
		Recall the definition of $f^*(\cdot)$ from~\eqref{eq:dualf} and $\partial f(\cdot)$ from~\cref{def:subg}. Then the following equivalence holds: 
		$$\langle x,y\rangle=f(x)+f^*(y) \quad \Longleftrightarrow \quad y\in\partial f(x)\quad \Longleftrightarrow\quad x\in \partial f^*(y).$$
	\end{lemma}
	\begin{lemma}[Bounding expected supremum of empirical process, see~\cite{arias2016estimation,Mason2012}]\label{lem:emproc}
		Let $f(\cdot)$ be a probability density supported on some subset of $\R^d$, and say $Z\sim f(\cdot)$. Let $\mathcal{G}$ be a class of uniformly bounded measurable functions from $\R^d\times (0,1]$ to $\R$, such that:
		$$\sup_{g(\cdot)\in\mathcal{G}}\sup_{h\in (0,1]} h^{-d}\E[g^2(Z,h)]<\infty,$$
		and such that the class
		$$\mathcal{G}_0:=\{x\mapsto g(x,h):\ g(\cdot)\in\mathcal{G},\ h\in (0,1)\}$$
		is point wise measurable and of VC-type (see~\cite[Section 2.6]{vaart2013} for relevant definitions of VC classes of sets/functions and point wise  measurability). Then there exists $b_0\in (0,1)$ such that if $Z_1,Z_2,\ldots $ is an i.i.d. sequence of observations from the probability density $f(\cdot)$, we have:
		$$\sup_{g(\cdot)\in\mathcal{G}}\sup_{\frac{\log{n}}{n}\leq h^d\leq b_0} \sqrt{\frac{n}{h^d\log{n}}}\Bigg|\frac{1}{n}\sum_{i=1}^n g(Z_i,h)-\E[g(Z,h)]\Bigg|=O_p(1).$$
	\end{lemma}	
		
\end{document}